\documentclass{amsart}

\usepackage{amssymb}
\usepackage{xypic}
\usepackage[dvipsnames]{color}
\usepackage{epsfig}

\newtheorem{theorem}[equation]{Theorem}
\newtheorem{lemma}[equation]{Lemma}
\newtheorem{proposition}[equation]{Proposition}
\newtheorem{prop}[equation]{Proposition} 
\newtheorem{corollary}[equation]{Corollary}
\theoremstyle{definition}
\newtheorem{algorithm}[equation]{Algorithm}
\newtheorem{assumption}[equation]{Assumption}

\newtheorem{convention}[equation]{Convention}
\newtheorem{definition}[equation]{Definition}
\newtheorem{example}[equation]{Example}
\newtheorem{method}[equation]{Method}
\newtheorem{notation}[equation]{Notation}
\newtheorem{remark}[equation]{Remark}

\theoremstyle{plain}
\numberwithin{equation}{section} 
\numberwithin{figure}{section}

\begin{document}

\title{Frobenius lifts and point counting for smooth curves}

\author{Amnon Besser}
\address{Amnon Besser, Mathematical Institute\\ 24--29 St Giles'\\ Oxford OX1 3LB\\ United Kingdom}

\author{Fran\c cois-Renaud Escriva}
\address{Faculteit der Exacte Wetenschappen\\Afdeling Wiskunde\\VU University Amsterdam\\De Boelelaan 1081a\\1081 HV Amsterdam\\The Netherlands}

\author{Rob de Jeu}
\address{Faculteit der Exacte Wetenschappen\\Afdeling Wiskunde\\VU University Amsterdam\\De Boelelaan 1081a\\1081 HV Amsterdam\\The Netherlands}

\begin{abstract}
We describe an algorithm to compute the zeta-function of a
proper, smooth curve over a finite field, when the curve is given
together with some auxiliary  data.
Our method is  based on computing
the matrix of the action of a semi-linear Frobenius on the first
cohomology group of the curve by means of Serre duality.
The cup product involved can be computed locally, after first
computing local expansions of a globally defined lift of Frobenius.
The resulting algorithm's complexity is softly cubic in the field degree,
which is also the case with Kedlaya's algorithm in the hyperelliptic case.
\end{abstract}

\subjclass[2010]{Primary: 14F30, 14G10, 14G15, 14Q50; secondary 14G22}

\keywords{curve over finite field, zeta function, rigid cohomology}

\maketitle

\def\({\left(}
\def\){\right)}

\def\iso{\simeq}
\def\leftiso{\buildrel{\simeq}\over{\leftarrow}}
\def\rightiso{\buildrel{\simeq}\over{\rightarrow}}

\def\ol{\overline}

\def\a{\alpha}
\def\b{\beta}
\def\cc{\gamma}
\def\d{\delta}
\def\e{\varepsilon}
\def\s{\sigma}
\def\o{\omega}

\def\D{\mathcal{D}}
\def\ee{\mathcal{E}}

\def\M{\mathbb{M}}

\def\A{\mathbb{A}}
\def\F{\mathbb{F}}
\def\H{\mathbb{H}}
\def\N{\mathbb{N}}
\def\Q{\mathbb{Q}}
\def\R{\mathbb{R}}
\def\Z{\mathbb{Z}}

\def\Zp{\Z_p}

\def\Cp{{\mathbb{C}_p}}

\def\Abar{\overline{A}}
\def\Adag{A^\dagger}
\def\phibar{\overline{\phi}}
\def\phirig{\phi_{\textup{cr}}^*}

\def\O{\mathcal{O}}
\def\OCq{\mathcal{O}_{C,q}}

\def\OmegaCF{\Omega_{\C_F}}

\def\dd{\textup{d}}
\def\wto{\widetilde\omega}

\def\tQ{\widetilde Q}
\def\tK{\widetilde K}
\def\tR{\widetilde R}

\def\X{\overline{X}}

\def\coker{\textup{coker}}
\def\Id{\textup{Id}}
\def\rig{\textup{rig}}
\def\dR{\textup{dR}}
\def\Jac{\textup{Jac}}
\def\MW{\textup{MW}}
\def\ord{\textup{ord}}
\def\Spec{\textup{Spec}}

\def\sump #1 #2 {\operatornamewithlimits{{\sum}^{\smash{\prime}}}_{#1}^{#2}}

\def\la{\langle}
\def\ra{\rangle}

\def\Rx{R\la \x \ra}
\def\Rxdag{R\la \x \ra ^\dagger}
\def\Rtab #1 #2 {R_{#1,#2}((t))}
\def\Rtstar{R_*((t))}
\def\Rt{R((t))}
\def\Rthat{\widehat{R((t))}}
\def\Stab #1 #2 #3 {S_{#1,#2}^{#3}((t))}

\def\Rts{R[[s,t]]}
\def\Ots{\O[[s,t]]}
\def\Rptab #1 #2 {R'_{#1,#2}((t))}

\def\ex{\xi}

\def\dr{{\textup{dr}}}
\def\Gal{\textup{Gal}}
\def\hdr{H_\dr}
\def\hcr{H_\textup{cr}}
\def\res{\operatorname{Res}}

\def\pair#1{\langle #1 \rangle}
\def\pairu#1{\pair{#1}_U}
\def\myq{q}

\def\pa{p'}

\def\fpb{{\overline{\F}_p}}

\def\Of #1 {\Omega^{1,f}_{#1}}

\def\ratpos{\Q_{>0}}
\def\x{\mathbf{x}}
\def\S{\mathbf{S}}
\def\sol{\mathbf{s}}
\def\soll{\tilde\sol}
\def\tsoll{\tilde\sol'}
\def\ts{\tilde{s}}
\def\tsp{\tilde{s}'}

\def\tN{\tilde N}

\def\Ia #1 {I_{#1}}
\def\Mlog{M_{\textup{log}}}

\def\ot{\tilde{O}}

\section{Introduction} \label{sec:intro}

Let $ p $ be a prime number and let  $ k $ be  a finite field of
cardinality $ q $ and characteristic $p$. An important problem of
algorithmic 
number theory is to count the number of points of a smooth (and
usually proper) variety $Y$ defined over $k$. By point
counting we mean, 
more precisely, the computation of the matrix of the $k$-linear
relative Frobenius map, acting on some \'etale or crystalline cohomology
group of $Y$. It is well-known that obtaining this matrix to a sufficiently
high precision allows an exact determination of its characteristic
polynomial as its coefficients satisfy the Weil bounds (see
Section~\ref{algorithm}).

The modern theory of point counting begins with the paper of
Schoof~\cite{Schoof85} for counting points on an elliptic curve $E$ by effectively
computing the action of Frobenius on the first \'etale cohomology group
of $E$. This direction of using \'etale cohomology is persued by various
other authors, still providing the best method when the field is prime
($p=q$) or close to being prime.

Other point counting methods, beginning with the work of Satoh~\cite{Sat00},
use crystalline cohomology. To describe these methods, let us fix some
more notation.

\begin{notation}\label{basicnot}
  Let $K$ be a finite Galois extension of the field $\Q_p$ of $p$-adic
  numbers, with ramification index $ e $, valuation ring $ R$,
  uniformizer $\pi$, and 
  residue field $R/\pi R$ isomorphic to $k$. We normalize the
  valuation on $ K $ by $v(p)=1$. We let $\sigma$ be an automorphism
  of $K$, and denote by $\ol{\sigma} $ the induced map on $ k $,
  which is given  by  $x\mapsto x^{\pa}$ with $\pa$ a positive power of $p$.
\end{notation}

For point counting one usually takes $e=1$ and $\pa=p$ but the theory
works in this generality,
and it will help us in future work
concerning syntomic regulators. We note that when $e=1$ the
automorphism $\sigma$ is uniquely determined by its reduction.

In counting methods based on crystalline cohomology, one computes an
effective representation for the crystalline cohomology group
$\hcr^i(Y/R) \otimes K$. This has a $\sigma$-semi-linear endomorphism
$\phirig$. It is obtained via functoriality of crystalline cohomology
from the relative Frobenius $\phibar$ on $Y$, i.e., the morphism
of schemes
\begin{equation*}
  \phibar: Y\to Y^{(\pa)}, \text{ where } Y^{(\pa)} =  Y\times_{k,\ol{\s}} k
\,,
\end{equation*}
obtained by raising to the $\pa$th power on the structure sheaf.
The sought after linear Frobenius is then obtained as a
(twisted) power of this.

When $Y$ can be lifted to characteristic $0$, and $e<p-1$, in
particular if $p$ is odd and $e=1$, 
crystalline cohomology
is the de Rham cohomology of the lift, and if the Frobenius
endomorphism can be lifted as well, then the endomorphism $\phi$ is
simply the action of the lift on this de Rham cohomology. This is
the case in, for example, Satoh's algorithm.

Finding lifts of Frobenius for a proper variety is rarely possible. An
alternative is to only lift Frobenius on an affine open piece. As
crystalline cohomology is infinite dimensional in this case, one has
to use a more refined cohomology theory, the Monsky-Washnitzer
cohomology~\cite{Mon-Was68}, which is a special case of Berthelot's
rigid cohomology~\cite{Ber97}. This cohomology theory associates to
the affine variety $\Spec(\Abar)$ the de Rham cohomology of
$\Adag_K = \Adag \otimes_R K$, where $\Adag$ is a ``weakly complete'' $R$-algebra whose
reduction modulo $\pi$ is $\Abar$. The action of $\phi$ is computed
from the action of a $\sigma$-semi-linear endomorphism of $\Adag$
reducing to the $p$-power map.

The use of Monsky-Washnitzer cohomology
in point counting algorithms was pioneered in the seminal paper
of Kedlaya~\cite{Ked01} on counting points on hyperelliptic
curves. Kedlaya's ideas can be extended to more general
curves~\cite{Ger03,Den-Ver06,CDV06} (see also the overview~\cite{Cha08}).

Kedlaya type algorithms generally consist of two main components.
\begin{enumerate}
\item An explicit lift of Frobenius to an endomorphism of $\Adag$,
 usually given in a straightforward manner.

\item A reduction algorithm that identifies a basis for
  $\hdr^i(\Adag_K)$ and shows how to explicitly write any $i$-form as
  a linear combination of basis elements plus an exact differential.
\end{enumerate}
Extending each of these steps from hyperelliptic curves to more
general curves in an efficient way proved to be a non-trivial task.

In this work we describe a point counting algorithm for curves under
the following fairly general assumptions:
we shall consider a proper, smooth curve 
$f:C\to\textup{Spec}(R)$
over $ R $ with geometrically irreducible fibres.
We shall denote by $C_K$ and $C_k$ its generic and special fibre
respectively.
Note that by Corollaire~7.4 of \cite[Expos\'e~III]{SGA1} we can lift any proper, smooth curve
over $ k $ to a smooth, proper curve $ C $ over $ W(k) $, necessarily
with geometrically irreducible generic fibre if $ C_k $ is geometrically
irreducible.

For the rest of this paper, we shall work with
the following situation and notations.
\begin{itemize}
\item
The genus of $C_K$ and $C_k$ is $g$.

\item
We are given a Zariski open affine $ X = \Spec(A) $ in $ C $ that
contains the generic point of $ C_k $, 
where $ A = R[x_1,\dots,x_n]/(f_2,\dots,f_n) $
with given generators $ f_2,\dots,f_n $ of the defining ideal;
moreover, the reduction $ \X = \Spec(\Abar) $ with
$ \Abar = A \otimes_R k = k[x_1,\dots,x_n]/(\ol{f_2},\dots,\ol{f_n})$
is a smooth complete intersection.
(See Assumption~\ref{Aassumption} for the terminology and
Remark~\ref{smoothlcirem} for the existence of such an~$ X $.)

\item 
We know $\omega_1,\ldots, \omega_{2g}$ in $ \Omega_{A/R}^1 $
that give a basis for the cohomology group $\hdr^1(C_K/K)$.

\item
We let $ \tK $ be a finite field extension of $ K $, with valuation ring
$ \widetilde R $, such that $ C_{\tR} \setminus X_{\tR} $
consists of the union of distinct sections $ Q_i : \Spec(\tR) \to C_{\tR} $,
which we tacitly identify with their images.
We denote the image of the closed point of $ \Spec(\tR) $ under
$ Q_i $ by $ q_i $. We do not assume that the $ q_i $ are distinct.
\end{itemize}

To perform point counting on these curves, we introduce three new
techniques. The first is a general explicit procedure for lifting
Frobenius in smooth, complete intersections situations
inspired by Section~2 of~\cite{Arabia01}. 
This method introduces a variable for each defining equation
$ f_j $, and uses
those to find a correction to the naive approximate lift given by raising
to the power $ p $. Using
this correction gives a map on $ R\la x_1,\dots,x_n \ra^\dagger $
that maps the defining ideal to itself and reduces to the desired
map on $ \Abar $.
We make this procedure explicit and provide 
estimates on the overconvergence of the resulting lift.

The second is a technique that avoids a 
(generally computationally expensive)
reduction algorithm by
replacing it with residue computations.
One observes that in order to know the matrix of $\phirig$
above,
it suffices to compute the cup products $\omega_i \cup \omega_j$
as well as the cup products $\phirig \omega_i \cup \omega_j$.
Our techniques reduce the computation of cup products to a
computation of residues of forms $(\int \omega_j) \phi \omega_i$ on
certain annuli, called ends, which are 
``at the boundary'' of the rigid space associated with $\Adag$.

Finally, essential for improving the performance of the algorithm, we
use a local lifting technique, which will compute the expansion of the
lifting of Frobenius locally near the boundary, instead of computing
it globally and restricting to the boundary.

Overall, the resulting algorithm for point counting is asymptotically
softly cubic in the field degree. This is the same complexity as
Kedlaya's algorithm~\cite{Ked01}, which is restricted to the case of
hyperelliptic curves, and the algorithm of Castryck, Denef and
Vercauteren for non-degenerate curves~\cite{CDV06}. The dependence on
the genus is somewhat worse for general curves but reduces in specific
situations. Finally, the dependence on $p$ is essentially linear. We
have not attempted an improvement in this direction in the style of~\cite{Har07}.

A (far from optimized) implementation of the above point counting alogrithm will be available within a few days.

The paper is organized as follows.
In Section~\ref{sec:cup-res} we explain how to obtain the matrix
of $ \phi $ from cup products on rigid analytic spaces.
In Section~\ref{global-frobenius} we discuss how to obtain the
desired lift $ \phi $ on $ \Adag $. Although we shall need it
only in the case of curves, we present the result and the estimates
on the coefficients involved for more general $ R $-algebras $ A $ 
that reduce to smooth complete intersections over $ k $.
Section~\ref{global-examples} makes the resulting maps and estimates
more explicit where $ X $ is an affine plane curve or a localization
of such a curve.  It also briefly discusses how to recover Kedlaya's
approach to hyperelliptic curves from our work.
Section~\ref{local-frobenius} discusses how to obtain the expansions
of the action of our lift at the ends, thus avoiding the computation
of the global lift constructed in Section~\ref{global-frobenius}.
Section~\ref{local-examples} returns to 
some of the examples discussed in Section~\ref{global-examples},
considering them from the points of view of leaving out only
one point, or obtaining a simpler lift $ \phi $ by localizing
more.
Section~\ref{sec:estimates} describes how to turn the estimates
of the preceeding theory into finite precision calculations that
still enable us to recover the zeta-function of $ C_k $, and
Section~\ref{algorithm} describes an algorithm  to do this, given
suitable input, and discusses its complexity.

Finally, we would like to thank
Bruno Chiarellotto,
Kiran Kedlaya,
Deepam Patel,
and
Jan Tuitman
for interesting and useful discussions.

\bigskip

Throughout the paper, we use the following notation.
\begin{notation}\label{not:1.2}
 We let $ R[[x_1,\dots,x_n]] $ denote
 the formal power series in $ x_1,\dots,x_n $ with coefficients in $
 R $, $ R\la x_1,\dots,x_n \ra $ the subring where the coefficients
 tend to 0 in $ R $, and $ R\la x_1,\dots,x_n \ra^\dagger $ the
 subring of $ R\la x_1,\dots,x_n \ra $ consisting of overconvergent
 power series.  We shall often use multi-index notation, writing $ \x
 $ for $ x_1,\dots,x_n $, and $ \x^I $ for $ x_1^{i_1} \dots
 x_n^{i_n} $ if $ I = (i_1,\dots,i_n) $. With $ |I| = i_1 + \dots +
 i_n $ we can then define $ \Rxdag $ as those $ \sum_I a_I \x^I $ in
 $ \Rx $ for which $ \cc $ and $ \d $ exist with $ \cc > 0 $ and
$ v(a_I) \ge \cc |I| + \d $ for all $ I $.
Equivalently, there exists a $ D $ in $ (\ratpos)^n $ and $ \d $
in $ \Q $ such that $ v(a_I) \ge D \cdot I + \d $ for all $ I $,
where $ D \cdot I  $ is the inner product.
\end{notation}

\section{Computing the matrix of Frobenius using cup products and residues} \label{sec:cup-res}

In this section we describe the
strategy for computing the matrix of Frobenius. For ease of
presentation we give a geometric description, based on
Coleman's work~\cite{Col-de88,Col89}. We then translate this into the more algebraic
language that will be used in the rest of the paper.

To rely directly on Coleman's work, it is convenient to first change
scalars to the field $\Cp$ of ``complex $p$-adic numbers''. Recall
that this is the completion of the algebraic closure of $\Q_p$. Its
residue field is the algebraic closure $\fpb$ of the finite field with
$p$ elements.

Consider $C_\Cp$ as a rigid analytic space over $\Cp$. Let $\myq_i$
be one of the $\fpb$ rational points in $ C_k \setminus \X $.
Let $\D_i \subset C_\Cp$ be
the rigid analytic subspace whose underlying set is the set of all
points whose reduction is $\myq_i$. This is called the residue disc of
$\myq_i$ by Coleman.

As $C$ is smooth, each of these $\D_i$ is
isomorphic to an open unit disc $\{z \in \Cp,\; |z|<1\}$. We choose
a parameter $t_i$ on
$\D_i$ realizing this isomorphism.
For each $0<r<1$ we let $U_r$ be the rigid subspace of $C_\Cp$ obtained
by removing the subsets $\{|t_i|\le r\}\subset \D_i$
(cf.~\cite[2.1]{Col-de88}). These $U_r$ are examples of ``wide open
spaces'' in Coleman's terminology.

We wish to compute the action of Frobenius on the cohomology of
$C_\Cp$. We assume we are given forms $\omega_1,\ldots,\omega_{2g}$ of
the second kind on
$C_\Cp$, whose cohomology classes form a basis of $\hdr^1(C_\Cp/\Cp)$, such
that all poles of the $\o_j$ are contained in the union of the $\D_i$.
By choosing a
sufficiently large $r_0$ we may assume that the $\o_j$ have no poles in
$U= U_{r_0}$.

We extend $\sigma$ to an automorphism of $\Cp$. We stress that this extension is not actually
used and is needed only so that we can formulate things over $\Cp$. We
will use a superscript $\sigma$ on the objects defined above to denote
the same object with structural morphism to $\Cp$ twisted by $\sigma$
(in the rigid analytic context this works better than twisted tensoring).
The following result also contains the definition of the ends
$ \ee_i $.

\begin{proposition} 
  There exist $r<1$, and with $U'= U_r$, a morphism $\phi:U' \to  U^\sigma$
whose reduction is the $\pa$-power map.  Furthermore, the morphism $\phi$
  has the property that $\phi(\ee_i^\prime) \subset \ee_i^\sigma$ with
  $\ee_i^\prime = U' \cap \D_i$ and $\ee_i^\sigma= U^\sigma\cap \D_i^\sigma$.
\end{proposition}

\begin{proof}
This is essentially~\cite[Theorem~2.2]{Col-de88} only in a semi-linear
version. The proof is the same.
\end{proof}

Recall that an annulus is a rigid analytic space isomorphic, via a
parameter $t$, to a space
of the form $\ee_r=\{r<|z|<1\}$. The space of rigid analytic functions
on such an annulus is
\begin{equation}
  \label{eq:annulus-fun}
  \left\{
\begin{aligned}
& \sum_{m\in \Z} a_m t^m\text{ with }a_m \text{ in } \Cp \text{ satisfying }
\\
&   \lim_{m\to -\infty} |a_m|  s^{-m} = 0 \text{ for all } s>r
     \text{ and}
\\
& \lim_{m\to \infty} |a_m|  s^{m} = 0 \text{ for all } s<1
\end{aligned}
\right\}
\,.
\end{equation}
The parameter $t_i$ restrict to
isomorphisms $\ee_i \to \ee_{r_0}$ and $\ee_i^\prime \to \ee_r$, so that both
$\ee_i$ and $\ee_i^\prime$ are annuli with parameter $t_i$.

\begin{definition}
  Let $\omega$ be a rigid analytic form on some annulus $\ee$
with parameter~$t$, and write
  \begin{equation*}
    \omega = \sum_{m\in \Z} a_m t^m \dd t\,.
  \end{equation*}
  Then we let the \emph{residue} of $ \omega $ on $\ee$ with respect to the
  parameter $t$ be $\res_\ee \omega = a_{-1}$.
\end{definition}

The set of all parameters on  an annulus $\ee$ breaks into two classes,
known as 
orientations~\cite[Lemma~2.1 and ensuing remarks]{Col89} such that the
residues with respect to
any two parameters are identical if they are in the same orientation
and differ by a
sign otherwise. An annulus with a choice of parameter in the same
orientation class is called an oriented annulus. The annuli $\ee_i$ are
oriented by the parameters $t_i$ and all parameters that are obtained
on the $\ee_i$ as restrictions of parameters on $\D_i$ give the same
orientation~\cite[Cor 3.7a]{Col89}. The choice of $t_i$ is therefore
irrelevant for the residue and we may denote it simply by
$\res_{\ee_i}$.

\begin{remark}\label{shrinkrm}
It is easy to see that for an annulus $\ee$ oriented by the parameter~$t$,
we have $\res_\ee = \res_{\ee'}$, where $\ee'$ is a subannulus defined
by the condition $ r' < |t| < 1 $.
\end{remark}

The analogue of the residue theorem holds \cite[Proposition~4.3]{Col89}.

\begin{theorem} \label{resthm}
  For a rigid analytic analytic form $\omega$ on $U$ we have
  $\sum_i \res_{\ee_i} \omega = 0$.
\end{theorem}

We recall the following basic result~\cite[Corollary~5.1]{Col98}.
\begin{theorem}\label{serre}
Let $ \o_1  $ and $ \o_2 $ be forms of the second kind on $C_\Cp$.
Then the cup product of their cohomology classes can be computed as
\begin{equation*}
    [\omega_1]\cup [\omega_2] = \sum_x \res_x \omega_2 \int \omega_1
\,,
\end{equation*}
  where the sum is over all points $x$ and the integral is a local
  integral with arbitrary constant term.
\end{theorem}

The key points to notice are that the integral makes sense since, with
respect to a local parameter $z$ at each point, there
is no term $az^{-1} \dd z$ to integrate, and that the constant of
integration does not
matter as in the residue computation it is going to multiply the
residue of $\omega_2$, which is $0$. If one of the forms is $\dd f$ for
a rational function $f$, then the residue theorem easily shows that
the right-hand side is indeed $0$.
\begin{definition}
  A rigid analytic form $\omega$ on $U$ will be called of the second
  kind if we have $\res_{\ee_i} \omega = 0 $ for every annulus
  $\ee_i$. If $\eta$ is another such form, the cup product pairing of
  $\omega$ and $\eta$ is defined by
  \begin{equation*}
    \pair{\omega,\eta} = \pairu{\omega,\eta} = \sum_i \res_{\ee_i} \eta
    \int \omega
  \,.
  \end{equation*}
\end{definition}

Just like in the algebraic setting of Theorem~\ref{serre}, it is clear
from the residue theorem, Theorem~\ref{resthm}, that the pairing is well-defined
and factors via $\hdr^1(U)$. It is further clear from
Remark~\ref{shrinkrm} that if $U'$ is a smaller wide open space as above
then 
\begin{equation*}
  \pair{\omega|_{U'},\eta|_{U'}}_{U'} = \pairu{\omega,\eta}
\,.
\end{equation*}

The usefulness of the pairing above for the computation of Frobenius
rests on the following result.
\begin{proposition}[{\cite[Proposition~4.5]{Col89}}]\label{cuppair}
  Let $\alpha_1,\alpha_2 \in \hdr^1(C_\Cp/\Cp)$ and let
  $\omega_1$ and $\omega_2$ be forms of the second kind on $U$ such
  that the class of $\omega_i$ in $\hdr^1(U/\Cp)$ is the restriction to
  $U$ of $\alpha_i$. Then $\alpha_1\cup \alpha_2 = \pairu{\omega_1,\omega_2}$.
\end{proposition}

The automorphism $\sigma$ acts on differential forms and cohomology
classes by sending them to the ``same'' forms and classes  on the
twisted objects. Note that as a vector space map, it is
$\sigma$-semi-linear, so that as expected, for example, it acts on differential
forms on an annulus $\ee$
with parameter $t$ by
\begin{equation*}
 \biggl(\sum_m a_m t^m \dd t \biggr)^\sigma
=
 \sum_m a_m^\sigma t^m \dd t
\,.
\end{equation*}

The cohomology group $\hdr^1(C_\Cp/\Cp)$ has a $\sigma$-semi-linear
endomorphism. Indeed, it is isomorphic to $\hcr^1(C/R)\otimes_R \Cp$
and the endomorphism is obtained by extending $\phirig$ on
$\hcr^1(C/R)$ $\sigma$-semi-linearly. We continue to denote this by $\phirig$.
To explicitly compute $\phirig$, we note that,
under restriction to $\hdr^1(U)$, it is compatible with the map
\begin{equation*}
  \hdr^1(U) \xrightarrow{\sigma} \hdr^1(U^\sigma)
  \xrightarrow{\phi^\ast} \hdr^1(U') \to \hdr^1(U)
\,,
\end{equation*}
where the last map is the inverse of the restriction map, which is an
isomorphism by~\cite[Theorem~4.2]{Col89}. This, and the compatibility of the
pairing with restrictions, immediately give the following.

\begin{corollary}
  Under the assumptions of Proposition~\ref{cuppair} we have
  \begin{equation*}
     \alpha_1 \cup (\phirig (\alpha_2)) = \pair{ \omega_1, \phi^\ast (\omega_2^\sigma)}_{U'}
\,.
  \end{equation*}
\end{corollary}

We can now describe our approach to computing the matrix of $ \phirig $.

\begin{method}
Assuming that one knows how to effectively compute the pairing
$\pairu{\,\cdot\, , \,\cdot\,}$, the above gives the following simple algorithm for
computing the matrix $ M $ of $ \phi^* $ with respect to the basis induced by
the $\{\omega_1,\ldots,\omega_{2g}\}$.

\begin{enumerate} \label{cupmethod}
\item Compute the cup product matrix $M_1$ with entries $\omega_i \cup
  \omega_j$ by using Theorem~\ref{serre}. 
\item Compute the cup product matrix $M_2=M_1 M$ with entries $
  \pair{ \omega_i,\phi^\ast \omega_j^\sigma}_{U'}$. 
\item Deduce the matrix $M$ of $\phirig$ as $ M_1^{-1}M_2 $.
\end{enumerate}
\end{method}

From this we can deduce the desired the zeta function.

\begin{method} \label{zeta-method}
Using Method~\ref{cupmethod}, compute the zeta function of $C_k$ as follows.
  \begin{enumerate}
  \item Compute the matrix $M$ of $\phirig$
    with a sufficiently high precision (see Section~\ref{algorithm}).
  \item Compute the matrix of the linear Frobenius as
    \begin{equation*}
      M' = \sigma^{l-1}(M)\times \sigma^{l-2}(M)\times \cdots \times
      \sigma(M)\times M
    \end{equation*}
    with $q=p^l$.
  \item Let $P_1(T)= \det(1-T M')$.
    Its coefficients are a priori in $\Z_p,$ but in fact are integers
    satisfying certain bounds deduced from the Weil bounds on the roots
    of $P_1$. Given a sufficiently high precision, $ P_1(T) $
    can therefore be determined precisely.
  \item Deduce the zeta fuction as $ Z(T) = \frac{P_1(T)}{(1-T)(1-qT)}$.
  \end{enumerate}
\end{method}

It remains to make concrete the computation of $
  \pair{ \omega,\phi^\ast \eta^\sigma}_{U'}$ for any two forms of the second
  kind on $U$. By our assumptions the parameters $t_i$ at the annuli
  can be chosen to be $\tK$-rational. The endomorphism $\phi$
  is induced by an endomorphism of dagger algebras
  ~\cite[2.2]{Col-de88}. The restriction of $\phi$ to $\ee_i^\prime$ is
  determined by the Laurent series expansion of $\phi^\ast t_i$, 
say  $f_i(t_i)$, with coefficients in $\tK$.
We then have
\begin{equation} \label{localcup}
  \pair{ \omega,\phi^\ast \eta^\sigma}_{U'}
=
  \sum_i \res_{\ee_i^\prime}  \phi^\ast \eta^\sigma \int \omega
\,,
\end{equation}
and given $f_i(t_i)$, each residue term is computed in
terms of Laurent series expansions
$ \omega = \sum_m a_m t_i^m \dd t_i $ and
$ \eta = \sum_m b_m t_i^m \dd t_i $.
It is the coefficient of $ t_i^{-1} $ in 
\begin{equation*}
 \biggl( \sum_m b_m^\sigma (f_i(t_i))^m \biggr) \biggl(\sump m {} \frac{a_m}{m+1}t_i^{m+1}\biggr)
\,,
\end{equation*}
where the prime denotes that we leave out the term with $ m=-1 $
in the sum (as $ a_{-1} = 0 $).
We shall describe more efficient methods for carrying out this
computation in later sections, but
at this point it is clear that it can be done in $\tK$.

\begin{remark} \label{dictionary}
Let us sketch the dictionary between this section and the rest of this
work, which is algebraic
rather than geometric. Rather than having a map of rigid spaces
$\phi:U' \to  U^\sigma$ we simply have a $\sigma$-semi-linear
endomorphism of the algebra $ \Adag$. We denote this by the same letter
$\phi$. This has the effect that the action on functions and
differential forms, which in this section is obtained by first
applying $\sigma$ to the coefficients and then applying $\phi^\ast$,
becomes in later sections simply the application of $\phi$ to the same
objects.
\end{remark}

\section{The global Frobenius} \label{global-frobenius}

Let $R$, $\pi$, $k$, $\sigma$ be as in Notation~\ref{basicnot}.
In this section we explain our strategy for computing a lift of
Frobenius on our dagger algebras over $R$, inspired by the work
of Arabia~\cite{Arabia01}. Even though we ultimately use this
only for curves, given the current limitation of our cup product
method for computing cohomology, the method applies, and we describe
it here, in greater generality for any $ A $ as in Assumption~\ref{Aassumption}.
By and large, this method was already developed in the master thesis of
F.-R.~Escriva. Later we discovered that another approach,
but with a less transparent presentation,
is contained in the unpublished PhD thesis of R.~Gerkmann~\cite{Ger03}.

Our goal in this section is to lift the $\pa$-power endomorphism $
\phibar $ of $ \Abar $ to a $\sigma$-linear endomorphism $\phi$ of
$\Adag= R \la x_1,\ldots,x_n \ra^\dagger /(f_{r+1},\ldots,f_n)$,
and obtain estimates on the coefficients of the $ \phi(x_i) $.
(See Remark~\ref{dictionary} for the relation with the notation
in Section~\ref{sec:cup-res}.)
We begin though, by explaining it in the
simplest possible case of one equation in two variables over $\Z_p$
and ignoring the issue of overconvergence.

Suppose then that we have $f(x,y)$ in $ \Z_p[x,y]$ such that the reduction
$\ol{f}(x,y)$ defines an non-singular curve in $\A^2$. Our goal is to
lift the Frobenius morphism $(x,y) \mapsto (x^p,y^p)$ to a morphism
$\phi$ of the affine curve defined by $f$, viewed as a rigid analytic variety.

Let $f_x$ and $f_y$ denote the partial derivatives of $f$ with respect
to the two variables. The non-singularity of $\ol{f}$ means that one
can find polynomials $\ol{P}_1$, $\ol{P}_2$ and $\ol{\Delta}$ in
$ \F_p[x,y]$
such that
\begin{equation*}
  \ol{P}_1 \ol{f}_x + \ol{P}_2 \ol{f}_y  = 1  +\ol{\Delta} \, \ol{f}
\,.
\end{equation*}
We arbitrarily lift $\ol{P}_1$, $\ol{P}_2$ and $\ol{\Delta}$ to polynomials
$P_1$, $P_2$ and $\Delta$ in $\Z_p[x,y]$, so that the congruence
\begin{equation}\label{congeq1}
  f_x P_1 + f_y P_2  \equiv 1+ \Delta f  
\end{equation}
holds modulo $ p $. We now seek our lift of Frobenius of the form
\begin{equation*}
  \phi(x,y) = (x^p,y^p) + s \times (P_1(x^p,y^p),P_2(x^p,y^p))
\end{equation*}
where $ s $ in $ p\Z_p \la x,y \ra$ is chosen to solve the equation
in the variable $S$,
\begin{equation}\label{newteq1}
  f\left( (x^p,y^p) + S \times (P_1(x^p,y^p),P_2(x^p,y^p))\right)
  - f(x,y)^p - f(x,y)^p \Delta(x^p,y^p) S = 0
\,.
\end{equation}
Clearly, if $s$ satisfies the above equation then $f(\phi(x,y))$ is
divisible by $f$ (even $f^p$), so that it indeed maps the curve defined
by $f$ to itself. Furthermore, since by assumption the coefficients of
$s$ are divisible by $p$, we see that $\phi(x,y) \equiv (x^p,y^p) $
modulo $ p$, so it is indeed a lift of Frobenius.

The equation~\eqref{newteq1} is an equation in one variable $S$ over
$\Z_p\la x,y \ra$, and $0$ is a solution modulo $p$. Its derivative with respect
to $S$ at $S=0$ is
\begin{equation*}
  f_x(x^p,y^p) P_1(x^p,y^p) + f_y(x^p,y^p) P_2(x^p,y^p)  - f(x^p,y^p)
  \Delta(x^p,y^p)
\,,
\end{equation*}
which, in light of~\eqref{congeq1}, reduces to $1$ modulo $p$. The
existence and uniqueness of the solution in $ p \Zp\la x, y \ra $
is thus guaranteed by Hensel's lemma, and it
can be recovered efficiently using Newton iterations starting from the
approximate solution $0$.

We shall consider the following very simple (and for point counting
obviously uninteresting) example at various
points in this paper in order to illustrate our estimates.

\begin{example} \label{runningexample1}
Consider $ f(x,y) = x^2 - y^2 -1 $ in $ \Zp[x,y] $ with $ p \ne 2 $.
Then $ \ol{2}^{-1} x \ol{f}_x + \ol{2}^{-1} y \ol{f}_y = \ol{1} + \ol{1} \cdot \ol{f} $
in $ \F_p[x,y] $. Now we write down
\begin{equation} \label{runexeq}
\begin{aligned}
   G(S)
& =
   f(x^p+2^{-1} x^p S, y^p + 2^{-1} y^p S) - f(x,y)^p - f(x,y)^p S 
\\
& =  
    4^{-1} (x^{2p} - y^{2p}) S^2 + (x^{2p} - y^{2p} -  f(x,y)^p ) S  - f(x^p,y^p) - f(x,y)^p
\end{aligned}
\end{equation}
in $ \Zp[x,y][S] $. We solve this for the unique solution $ S = s $
in $ p \Zp\la x, y \ra $. Then $ \phi(x,y) = (x^p+2^{-1}x^p s, y^p+2^{-1}y^p s) $
induces an endomorphism of $ \Zp\la x,y \ra $ that descends
to an endomorphism of $ \Zp\la x,y \ra  /(f(x,y)) $ because
it maps the ideal $ (f(x,y)) $ to itself by construction, and it reduces
to the Frobenius map $ \phibar(x,y) = (x^p,y^p) $ modulo $ p $.
\end{example}

We now describe the general case, still ignoring
overconvergence. 

Recall the shorthand $ R[\x] $ of Notation~\ref{not:1.2}.
We shall also write $\M^{a,b} $ and $\M^{a}$ for $a\times b$ and
$a\times a$ matrices, and if $ f_{r+1}, \dots, f_n $ in $ R[\x] $
are given, then we let
$\Jac_f $ in  $ \M^{n-r,n}(R[\x])$ be the resulting Jacobian matrix.
We shall lift  the $\pa$-power endomorphism $ \phibar $ of $ \Abar $ to a $\sigma$-linear endomorphism $\phi$ of
$\Rx/(f_{r+1},\ldots,f_n) $ for the following $ A $.
In particular, by Remark~\ref{smoothlcirem} below,
this will apply to a suitable Zariski open part $ X $ of $ C/R $.

\begin{assumption} \label{Aassumption}
In $ R[\x] $, for $ 0 \le r \le n-1 $,
we are given $ f_{r+1}, \dots, f_n $,
such that 
\begin{equation} \label{Apresentation}
 A = R[\x] / (f_{r+1}, \dots, f_n) 
\,.
\end{equation}
If  $ \ol{\Jac_f} $ is the reduction modulo $ \pi $ of $ \Jac_f $, then
the unit ideal in $\Abar$ is generated by the determinants
of the $(n-r) \times (n-r)$ minors of $ \ol{\Jac_f} $.
\end{assumption}

Under this assumption, Arabia shows in the proof of~\cite[Th\'eor\`eme~2.1.2]{Arabia01}
that there exist matrices 
\begin{align}
\label{defpdel}
  P \in \M^{n,n-r}(R[\x])\;&, \quad \Delta^{r+1},\ldots,\Delta^n \in  \M^{n-r}(R[\x])
\\
 \intertext{such that}
\label{congeq2}
    \Jac_f \times P &\equiv \Id_{n-r}+\sum_{j=r+1}^n f_j \Delta^j \textup{ modulo }\pi
\,.
\end{align}

Let $\psi$ be the $\s$-linear endomorphism of $ \Rx $ that sends each
$ x_i $ to $ x_i^{\pa} $, so that it maps
an element $g(\x)$ to $g^\s(\psi(\x))$, where the superscript
$ \s $ means we apply $ \s $ to the coefficients.
 We shall look for a $ \s $-linear
$\phi$, defined by its action on the column vector of variables $\x$ as
\begin{equation} \label{eq:phiact}
  \phi(\x) = \psi(\x)+\psi(P)\sol
\,,
\end{equation}
where $\sol$ is then a column vector in $\pi \Rx^{n-r}$.
We want $ \sol $ to satisfy $G(\sol)=0$, where
the column vector $ G(\S) =  (G_{r+1}(\S),\dots, G_n(\S)) $
with entries in $ \Rx [\S] $, is given by
\begin{equation} \label{globalG}
  G(\S) = f^\s (\psi(x)+\psi(P)\S) - f^{\pa} -
 \sum_{j=r+1}^n f_j^{\pa} \psi(\Delta^j) \S  
\,,
\end{equation}
for $ f^\s $ the vector $ (f_{r+1}^\s,\dots,f_n^\s) $,
$ f^{\pa} $ the vector $ (f_{r+1}^{\pa},\dots,f_n^{\pa}) $,
and $ \S $ the vector $ (S_{r+1},\dots,S_n) $.

In a way similar to the case of one equation in two variables discussed
before, one finds that
\begin{itemize}
\item
$G(0) \equiv 0 $ modulo $ \pi $;

\item
$ \Jac_G(0)  = \Jac_{f^\s}(\psi(x)) \times \psi(P) - \sum_{j=r+1}^n f_j^{\pa} \psi(\Delta^j) $,
hence $\Jac_G(0) \equiv\Id_{n-r} $  modulo $ \pi $.
\end{itemize}
Therefore the equation may be solved uniquely
for $ \sol $ in $ \pi \Rx^{n-r} $ by Hensel's Lemma,  and this
can be done
effectively using Newton iteration. 
It is now clear that $\phi$ is $ \s $-linear, reduces to the
$ \pa $-power map $ \phibar : k[\x] \to k[\x] $,
and maps
the ideal $ (f_{r+1},\dots,f_n) $ into itself. Overall, we obtained the
following result.

\begin{theorem}\label{globallift0}
  Let $f_{r+1}, \dots, f_n $ in $ R[x_1,\ldots, x_n]$
  with $ 0 \le r \le n-1 $ be given, and suppose
  that $A=R[x_1,\ldots,x_n]/(f_{r+1},\ldots,f_n)$ satisfies
  Assumption~\ref{Aassumption}.
  Fix $ P $ and $ \Delta^{r+1}, \dots, \Delta^n $
  as in~\eqref{defpdel} and~\eqref{congeq2},
  and let $G(\S)$ be given by \eqref{globalG}.
  Then there exists a unique lift of the
  $\pa$-power map on $\Abar$ to a $\sigma$-semi-linear endomorphism of
  $\Rx/(f_{r+1},\ldots,f_n)$ of the form~\eqref{eq:phiact}
  with $\sol $ in $ \pi \Rx^{n-r}$ satisfying $G(\sol)=0$.
\end{theorem}

To be able to effectively use this lift of Frobenius, we need to know
that it preserves overconvergence, and know explicit bounds on the rate of
convergence. To this end, we first need explicit bounds on the rate of
convergence obtained in Hensel's Lemma. These are provided by
Lemma~\ref{polytopelemma} below.
In order to describe them, it will be convenient
to introduce some notation.

\begin{notation}
For non-zero $ D $ in $ \Q_{\ge 0}^n $, let
$ V_D = \{ v \textup{ in } (\R_{\ge0})^n \textup { with } D \cdot v \le 1 \} $.
\end{notation}

In the process of obtaining our estimates in
Lemma~\ref{polytopelemma} and similar results in
Section~\ref{local-frobenius}, we shall introduce suitably ramified
extensions. In order to avoid interrupting the flow of the argument,
we impose the following.

\begin{convention} \label{ratconv}
If $ \a $ is in $ \ratpos $, then $ \pi^\a $ means that we extend
the ring $ R $ to the valuation ring $ R' $ in a finite extension
$ K' $ of $ K $ for which $ \a $ is attained as a valuation.
In other words, such that $ \pi^\a $ can be interpreted as a
integer power of a uniformizer of~$ R' $.
\end{convention}

\begin{lemma} \label{polytopelemma}
Let $G_{r+1}(\S),\dots,G_n(\S)$ in $R[\x][\S]$ be of maximal total
degree $ N > 0 $ in $ \S $. For $ l=0,\dots,N $, let $ G_{j,l}(\S) $
consist of the terms of $ G_j(\S) $ that are homogeneous in $ \S $
of degree~$ l $.
Assume that the $G_j(0)$ are in $\pi R[\x]$ and that the determinant
of $ \Jac_G(0) $ is in $ R^* + \pi R[\x] $.
Then there is a unique solution $ \sol $ in $ (\pi \Rx )^{n-r} $
of $G_{r+1}(\S) = \dots = G_n(\S) = 0 $; in fact, it
lies in $ (\pi \Rxdag )^{n-r} $.

Moreover, $ \sol = \sum_I a_I x^I $ where for each coordinate
$ a_{I,j} $ we have the following estimate, independent of $ j $.
If $ D $ in $ (\Q_{\ge0})^n $, and  $a, b$ in $ \ratpos $ are such that
the Newton polytope of $ G_{j,l}(\S) $ is contained in $(al+b) V_D$
for $ l=0,\dots,N $, then for each $ I $ we have
\begin{equation*}
 v(a_{I,j}) \ge \frac{D \cdot I +2a+b}{2e(a+b)}
\,.
\end{equation*} 
\end{lemma}

\begin{proof}
View $ G = (G_{r+1},\dots,G_n) $ as column vector.
Then the existence and uniqueness of $\sol$ in $ (\pi \Rx)^{n-r} $
are obtained from Hensel's lemma, by starting with $ z_0 = (0,\dots,0) $
as approximate solution of the vector equation $ G(\S) \equiv 0 $ modulo $ \pi $,
and applying Newton iteration
$ z_{i+1} = z_i - \Jac_G(z_i)^{-1} G(z_i) $ for $ i \ge 1 $.

For the estimate, let $\e $ be in $ \ratpos $ and define
$   \mu =\frac{1}{2(a+b)+\e} D \textup{ in } \Q_{\ge 0}^n \textup{ and }\nu=\frac{2(a+b)- b}{2(a+b)+\e} $
in $  \ratpos $.
Then 
\begin{equation} \label{munucond}
\begin{aligned}
 & \mu \cdot I < 1-\nu \textup{ for all } I \textup{ in } bV_D;
\cr
& \mu \cdot I < \frac{1}{2} \textup{ for all } I \textup{ in } (a+b)V_D;
\cr
& \mu \cdot I \le (l-1)\nu \textup{ for all } I \textup{ in } (al+b)V_D \textup{ with } l=2,\ldots,N
\,.
\end{aligned}
\end{equation}

Using Convention~\ref{ratconv} above,
we apply to each $ \pi^{-\nu}G_j(\S)$ the substitutions
$S_j\leftarrow\pi^{\nu}S_j'$ and $x_i \leftarrow \pi^{-\mu_i} x_i'$.
We shall abbreviate the latter to $ \x \leftarrow \pi^{-\mu} \x $.
In order to describe the result we abuse notation and write $ G_{j,l}(x,\S) $
for $ G_{j,l}(\S) $.  Then we obtain 
\begin{equation*}
  G_j(\S')
 = \sum_{l=0}^N \pi^{\nu(l-1)} G_{j,l}( \pi^{-\mu}\x' ,\S' ) 
\end{equation*}
in $ K'[\x'][\S'] $.
Applying the first inequality in~\eqref{munucond} to all $G_{j,0}$, the
second to all $G_{j,1}$, and the third to all
$G_{j,l}$ for $l\geq 2$, one sees that each $G_j(\S')$ is in
$ R'[\x'][\S'] $.  Moreover, each $G_{j,0}$ is in $\pi'R'[\x']$,
and the determinant of $ \Jac_G(0) $ is in  $R'^* + \pi' R'[\x'] $.
By Hensel's lemma there exists a unique solution  
$ \sum_{I\ge 0} b_I \x'^I $ in $(\pi'R'\la \x'\ra)^{n-r}$ of
$  G_{r+1}(\S') = \dots =  G_n(\S') = 0 $.
Thus
$ \sum_I \pi^{\mu \cdot I+\nu} b_I \x^I $
and $ \sol $ are two solutions of $G_{r+1}(\S)=\cdots = G_n(\S)=0$
in $(\pi'R'\la x\ra)^{n-r}$,
but by Hensel's lemma in $ R'\la x \ra $ there is only one such solution.
So for each $ I $ we find
\begin{equation*}
  v(a_I) > \frac{1}{e}\( \mu \cdot I + \nu \) = \frac{ D \cdot I + (2a+b)}{e(2(a+b)+\e)}
\,.
\end{equation*}
Letting $ \e $ go to~0 we find 
$ v(a_I) \ge \frac{D \cdot I + (2a+b)}{2e(a+b)} $.
Using $ (1,\dots,1) $, $ a $ sufficiently large and $ b=0 $
shows that $ \sol $ is in $ (\pi \Rxdag )^{n-r} $.
\end{proof}

We can now lift the endomorphism $ \phibar $ of $ \Abar $ to an endomorphism of $ \Adag $.

\begin{theorem} \label{globallift}
Let $\phi$ be the lift of Frobenius constructed in
Theorem~\ref{globallift0}. Then $\phi$ preserves $\Adag$. Furthermore,
the following estimate holds for the coefficients $ b_I = b_{i,I} $ in
each $ \phi(x_i) = x_i^{\pa} + \sum_I b_{i,I} x^I $, 
and is independent of $ i $.
Let $\Gamma$ be the intersection of all $ V_D $
that contain the Newton polytopes of all $f_j$.
Fix $d$ in $ \ratpos $ such that the Newton polytopes
of all coefficients in the matrices $ P $ and
$\Delta^{r+1},\dots,\Delta^n$ in~\eqref{defpdel}  are 
included in $d\Gamma$.
If $ c $ is in $ \ratpos $, then
$v(b_I)>\frac{c}{2e(d+1)\pa}+\frac{1}{2e}$
whenever $I$ is not in~$c\Gamma$.
\end{theorem}

\begin{proof}
Let $G(\S)$ be as in~\eqref{globalG}.
We recall that the lift $\phi$ is given by the formula
\begin{equation*}
  \phi( g(\x) ) = g^\s \(\psi(\x) + \psi(P) \sol \) 
\end{equation*}	
using the unique
$\sol$ in $ (\pi \Rx)^{n-r}$ with $ G(\sol) = 0 $.

We now prove the estimate of the coefficients,
which will also show that $\Adag$ is preserved.
Note that $ \Gamma $ automatically contains the Newton polytopes
of all (higher) partial derivatives of all $ f_j $.
Then one checks easily that
the Newton polytopes of the entries of the homogeneous part of degree $ l $
in $ \S $ of $ G(\S) $ are contained in $(dl+1){\pa}\Gamma$.

Write $s_j = \sum_I a_{j,I} x^I $. If $ I $ is not in $ c \Gamma $,
then $ D \cdot I > c $ for some $ D $ in $ (\Q_{\ge0})^n $ with
$ \Gamma \subseteq V_D $.
Then by Lemma~\ref{polytopelemma}, with $ a=d \pa $ and $ b=\pa $,
we have
\begin{equation} \label{overconvergenceestimate}
     v(a_{j,I}) \ge \frac{D \cdot I +\pa(2d+1)}{2e(d+1)\pa}  >  \frac{c +\pa(2d+1)}{2e(d+1)\pa}
\,.
\end{equation}
Fix $ i $, and write
$\psi(P_{i,j}) = \sum_{K\in \pa d\Gamma} d_{j,K} x^K$.
Then $ \phi(x_i) = x_i^{\pa} + \sum_{j=r+1}^n \psi(P_{i,j}) s_j = x_i^{\pa} + \sum_I b_I x^I $
with
\begin{equation*}
b_I = \sum_{j=r+1}^{n}{\sum_{K+L=I}{d_{j,K}a_{j,L}}}
   \,.
\end{equation*}
Take $ c $ in $ \ratpos $ and assume $ I $ is not in $ c \Gamma $.
If $ c \le \pa d $, then $v(b_I)\ge \frac{1}{e}>\frac{c}{2e(d+1)\pa}+\frac{1}{2e}$
because $ a_{j,L} $ is in $ \pi R $.
If $c > \pa d $, then in $ I = K+L $ we have that $L$ is not in
$(c-d\pa)\Gamma$.  
By~\eqref{overconvergenceestimate} we then have
$  v(a_{j,L}) > \frac{(c-\pa d) +\pa(2d+1)}{2e(d+1)\pa}  = \frac{c}{2e(d+1)\pa}+\frac{1}{2e} $.
Because all $ d_{j,K} $ are in $ R $, our estimate has been proved.
\end{proof}
 
\begin{remark}\label{liftcurve}
In Theorem~\ref{globallift} one can sometimes
prescribe that $ \phi(x_i) = x_i^{\pa} $
for some $ i $. 
If $ A $ is given by a presentation
\begin{equation*}
   0\rightarrow (f_{r+1},\dots,f_n)\rightarrow R[\x]\rightarrow A\rightarrow 0
\end{equation*}
and $ a $ is a positive integer with $ a \le r $,
let us denote by $ J_a $ the matrix consisting of the last $ n - a $ columns of
$ \Jac_{(f)} $.
If the $ (n-r) $-minors of $ \ol{J_a} $ generate
the unit ideal of $ \Abar $, then one can compute a lift $ \phi $ of
$ \phibar $ with $ \phi(x_i) = x_i^{\pa} $
for $ i=1,\dots,a $. 
Namely, we can apply the result of Arabia \eqref{congeq2} with
our $ R $ replaced with $ R[x_1,\dots,x_a] $, and $ R[x_1,\dots,x_n] $
with $ R[x_1,\dots,x_a][x_{a+1},\dots,x_n] $. We then obtain matrices
$ P_a $ in $ \M^{n-a,n-r}((R[x_1,\dots,x_n])) $
and $ \Delta^{r+1}, \dots, \Delta^n $ in $ \M^{n-r}((R[x_1,\dots,x_n])) $
with $ J_a \times P_a  \equiv \Id_{n-r} + \sum_{j=r+1}^n f_j \Delta^j $
modulo $ \pi $.  This means we have satisfied~\eqref{congeq2}
with a matrix $ P $ for which the first $ a $ rows are identically~0,
and~\eqref{globalG} now becomes
\begin{alignat*}{1} 
  G(\S)
& =
  f^\s\(\psi(x_1),\dots,\psi(x_a),\psi(x_{a+1})+\psi(P_{a+1,*})\S,\ldots,\psi(x_n)+\psi(P_{n,*})\S\)
\\
& \phantom{\,=\,} -f_j^{\pa} -\sum_{j=r+1}^{n} f_j^{\pa}\psi(\Delta^j)\S
\,.
\end{alignat*}
\end{remark}

We return to our running example, Example~\ref{runningexample1},
as an illustration of the estimates in Theorem~\ref{globallift}.

\begin{example} \label{runningexample2}
We apply the estimates of Theorem~\ref{globallift}
to the data in Example~\ref{runningexample1}, so that we have
$ m=2 $, $ r = 1 $,
$ P = \begin{pmatrix} 2^{-1} x \\ 2 ^{-1} y \end{pmatrix} $,
$ f_2(x,y) = x^2 - y^2 - 1 $, 
$ \Delta^2 = \begin{pmatrix} 1 \end{pmatrix} $,
$ \Gamma = V_D $ for $ D = (\frac12,\frac12) $,
and $ d = \frac12 $.
Then $ (i,j) $ is not in $ c \Gamma $ if and only if
$ i +j > 2 c $.  Therefore,
if $ \phi(x) = x^p + \sum_{i,j \ge 0} b_{i,j} x^i y^j $,
then $ i+j > 2c $ implies $ v(b_{i,j}) > \frac{c}{3p} + \frac12 $.
Equivalently, $ v(b_{i,j}) \ge \frac{i+j}{6p} + \frac12 $ for all $ (i,j) $.
The same estimates holds for the coefficients $ b_{i,j}' $ in
$ \phi(y) = y^p + \sum_{i,j} b_{i,j}' x^i y^j $.
\end{example}

We conclude this section by showing that our theorems apply
to suitable open parts of smooth, Noetherian schemes over $ R $.

\begin{remark} \label{smoothlcirem}
Suppose that $ Y $ is a smooth, Noetherian scheme $ Y $ over $ R $ of
relative dimension $ r \ge 0 $. Then
there exists a Zariski open affine part that is of
the form~\eqref{Apresentation},
and such that the unit ideal of $ A $ is generated by the determinants
of the $ (n-r) \times (n-r) $ minors of $ \Jac_f $. Moreover,
there exist matrices $ P $ in $ \M^{n,n-r}(R[\x]) $ and $ \Delta^{r+1},\ldots,\Delta^n$
in $ \M^{n-r}(R[\x]) $ such that
\begin{equation*}
    \Jac_f \times P = \Id_{n-r}+\sum_{j=r+1}^n f_j \Delta^j 
\end{equation*}
in $ \M^{n-r} (R[\x]) $.

Namely, let $ y $ be the generic point of the special fibre $ Y_k $. 
Since $ Y $ is smooth over $ R $ of relative dimension $ r $, there exists
an open neighbourhood $ U $ of $ y $ in $ Y $ and an immersion $ j $ of $ U $ 
into an affine space $ \A_R^m $, such that, locally around
$ z = j(y) $, the ideal sheaf defining $ j(U) $ in some open of 
$ \A_R^m $ is generated by $ m - r $ sections $ g_{r+1},\dots, g_m $. 
Furthermore, the differentials $ \dd g_{r+1}(z),\dots, \dd g_m(z) $ are linearly 
independant in $ \Omega^1_{\A_R^m/R}\otimes k(z) $. Note that every open 
$ V\subset U $ containing $ y $ also has this property.

According to \cite[Proposition~(17.2.5)]{EGA4}, after localizing
more if necessary, 
there exists such an open affine neighbourhood of $ y $ on which
the conormal exact sequence splits. 
We may assume it is given by an algebra $ A = R[\x]/J $, for an ideal 
$ J = (f_{r+1},\dots, f_n) $ of $ R[\x] $.
Then the morphism $ \d $ in the exact sequence of $ A $-modules
\begin{equation*}
	\xymatrix{
	J/J^2 \ar[r]^-\d
	& \Omega^1_{R[\x]/R}\otimes_{R[\x]} A \ar[r]
	& \Omega^1_{A/R} \ar[r]
	& 0
	}
\end{equation*}
is injective and admits a retraction.
Therefore, $ A\cdot\dd f_{r+1} \oplus \dots \oplus A\cdot\dd f_n \cong A^{n-r} $ 
is a direct summand of $ \Omega^1_{R[\x]/R}\otimes_{R[\x]} A\cong A^n $, and
there exists a right inverse $ P $ in $ \M^{n,n-r}(A) $ of $ \Jac_f $.
\end{remark}

\section{Examples of the global Frobenius} \label{global-examples}

In this section we make the construction of $ \phi $ in
Theorem~\ref{globallift}
more explicit in the case of plane curves and their localisations.
Note that then $ r=1 $ and $ n=2 $ or 3.

\begin{example} \label{planeexample}
Let us first treat the case of a smooth curve in $ R[x,y] $,
defined by $ f(x,y) $, with the current notation.
There exist $ P_1 $, $ P_2 $ and $ \Delta $ in
$ R[x,y] $ such that
$ \ol{P_1} \, \ol{f_x} + \ol{P_2} \, \ol{f_y} = 1 + \ol{\Delta} \cdot \ol{f} $
in $ k[x,y] $, and~\eqref{globalG} becomes
\begin{alignat*}{1}
 G(S)
& =
f^\s \bigl( x^{\pa} +  \psi(P_1) S, y^{\pa} + \psi(P_2) S \bigr) -
f^{\pa} - f^{\pa} \psi(\Delta) S 
\\
& =
     \psi(f) - f^{\pa}
   + \(f^\s_{x}\bigl(x^{\pa}, y^{\pa}\bigr) \psi(P_1) 
                   + f^\s_{y}\bigl(x^{\pa}, y^{\pa}\bigr)\psi(P_2) - f^{\pa} \psi(\Delta) \) S
\\
& \phantom{\,=\,}   + (\cdots) S^2 + \cdots
\,.
\end{alignat*}
With $ s $ the unique solution in $ \pi R\la x, y \ra ^\dagger $
of $ G(S) = 0 $,
the map from $\Adag = R\la x \ra^\dagger/(f)$ to itself is given
by mapping the class of $x$ to that of $x^{\pa} + P_1^\s(x^{\pa} , y^{\pa}) s $,
and the class of $y$ to that of $y^{\pa} + P_2^\s(x^{\pa} , y^{\pa}) s $.

Moreover, let $\Gamma$ be the intersection of all $ V_D $ that
contain the polytope of $f$, and let $d$ be a positive
rational number such that the Newton polytopes of $P_1$, $P_2$ and $\Delta$ are
included in $d\Gamma$.
Then $ \phi(x) = x^{\pa} + \sum_{i,j\ge 0} b_{i,j} x^i y^j $ where, for each
positive rational number $c$, 
we have $v(b_{i,j})> \frac{c}{2e(d+1)\pa} + \frac{1}{2e} $
whenever $(i,j)$ is not in $ c\Gamma$.
The same estimates apply to the coefficients in $ \phi(y) - y^{\pa} $.
\end{example}

\begin{example} \label{ellipticexample}
Let us treat an explicit case of Example~\ref{planeexample}.
Consider the elliptic curve over $\Zp$ with $ p \ne 2,3 $
defined by $f(x,y) = y^2-x^3-1$. Then we even 
have $ \frac13x f_x + \frac12 y f_y = 1 + f$ in $ \Zp[x,y] $.
Let us take $ \pa = p $.
Noting that $ \s $ is the identity here,
we have to find the unique solution $ s $ in $ p\Zp\la x,y \ra^\dagger $
of $ G(S) = 0 $, with $ G(S) $ the polynomial
\begin{alignat*}{1}
 G(S)
& =
 f\left( x^p + \tfrac13 x^p S, y^p + \tfrac12 y^p S \right) - f^p-f^p S 
\\
& =
 -\tfrac{1}{27} x^{3p} S^3 + \left(\tfrac14 y^{2p} - \tfrac13 x^{3p}\right)
   S^2 + (1+f(x^p,y^p) - f^p) S + f(x^p,y^p)-f^p
\,.
\end{alignat*}
The map of $ \Zp \la x,y\ra^\dagger/(y^2-x^3-1) $ to
itself is then given by mapping the class of $ x $  
to that of $ x^p + \frac13 x^p s $,
and the class of $ y $ to that of $ y^p + \frac12 y^p s $.
The polytope $ \Gamma $ has vertices $ (0,0) $, $ (3,0) $ and $ (0,2) $,
and equals $ V_D $ with $ D = \frac16 (2,3) $.
We can take $ d= \frac12 $, so
that $ \phi(x) = x^p + \sum_{i,j\ge0} b_{i,j} x^i y^j $ with
$ v(b_{i,j}) > \frac{c}{3 p} + \frac12 $ whenever $ 2i+3j >6c $.
In fact, $ (i,j) $ is not in $ c\Gamma $ if and only if $ 2i+3j > 6c $,
so choosing $  c= \frac{2i+3j}{6} - \e $ with $ \e>0 $ and letting
$ \e $ go to 0 we find that
$ v(b_{i,j}) \ge \frac{2i+3j}{18p} + \frac12 $.
The same estimates apply to the coefficients in $ \phi(y) - y^p $.
\end{example}

\begin{example} \label{raisex1example}
Let us consider an irreducible affine curve defined by $ f(x,y) $ with $ f_y $
not identically 0 modulo~$ \pi $.
Let $ A = R[x,y,z]/(f,zf_{y} - 1) $.
Notice that $ A $
satisfies the assumption of Remark~\ref{liftcurve} with $ n = 3 $ and $ r = k = 1 $.
Let $ h $ denote $\frac{\partial^2f}{\partial y^2}$. Then we
have
$ \Jac 
=
\begin{pmatrix}
  f_{x}     & f_{y} & 0       \\
  z f_{x,y} & z h   & f_{y} \\
\end{pmatrix}
$,
so we can take
$ P = \begin{pmatrix}
    0    &  0    \\
    z    &  0    \\
  -z^3 h &  z  \\
\end{pmatrix}
$,
$ \Delta_2 = 0 $
and
$ \Delta_3 = 
\begin{pmatrix}
  1     & 0  \\
 -z^2 h & 1 \\
\end{pmatrix}
$.
Thus, we have to find the unique solution 
$ (s_2 , s_3) $ in $(\pi R\la x,y,z \ra^\dagger)^2$ of the equations
\begin{alignat*}{1}
   G_2(S_2, S_3) 
   & = f^\s\(\psi(\x) + \psi(P)\S\) - f^{\pa} - \(zf_y-1\)^{\pa} S_2\\
& =
   f^\s\(\psi(x), \psi(y)+\psi(z)S_2\)
   - f^{\pa} 
   - \(z f_y -1 \)^{\pa} S_2
\\
\intertext{and}
   G_3(S_2, S_3)
& =  
   f_3^\s\(\psi(\x) + \psi(P)\S\)
   - \(zf_{y}-1\)^{\pa} \( 1 - \psi(z^2 h)S_2 + S_3 \) 
\end{alignat*}
with $ f_3(x,y,z) = z f_y -1 $.
Note that the first term in $ G_3(S_2, S_3) $ is then given explicitly
as
$ (\psi(z) - \psi(z^3h) S_2 + \psi(z) S_3) f_y^\s(\psi(x), \psi(y)+\psi(z) S_2)-1 $.

Observe that $ G_2 $ is a polynomial in $S_2$ only, and that
by Lemma~\ref{polytopelemma} there exists
a unique solution $ s_2 $ in
$ \pi R\la x,y,z\ra ^\dagger $ of $ G_2(S_2) = 0 $.
Then the map $ \phi $ from $ R \la x, y, z \ra ^\dagger $
to itself maps $ x $ to
$x^{\pa}$, and $ y $ to $ y^{\pa} + z^{\pa}s_2 $.

Let us notice that all the coefficients of $ P $ and $ \Delta_3 $ lie in 
$ 3 \Gamma $, where $ \Gamma $ is the intersection of all $ V_D $ that contain
the Newton polytopes of $ f $ and of $ zf_y - 1 $.
Thus, if we write 
$ \phi(y) = y^{\pa} + \sum_I b_I x^I $, then for every $ I $ that is not in 
$ c\Gamma $ with $ c $  in $ \ratpos $ fixed, we have the estimate 
$ v(b_I) > \frac{c}{8e\pa} + \frac{1}{2e}$. 

In order to determine $ \phi(z) $ in $ R\la x, y, z \ra^\dagger $
one would have to solve the
equations, but the class of $ \phi(z) $ in $ \Adag $ is determined
by $ \phi(z) \phi(f_y)=1$. Note that 
\begin{equation*}
  \phi(f_y) = f_y^\s\bigl(x^{\pa}, y^{\pa}+z^{\pa}s_2\bigr)
=
  \psi(f_y) + F_1
=
  f_y^{\pa} - F_2
=
  z^{-\pa} - F_2
\,,
\end{equation*}
where $ F_1 $ and $ F_2 $ are in $ \pi R\la x,y,z\ra^\dagger $.
Therefore the class of $ \phi(z) $ equals that of
$ z^{\pa} \bigl(1+\sum_{m=1}^\infty (z^{\pa}F_2)^m\bigr)$.
\end{example}

\begin{example} \label{kedlayaexample}
Let us apply Example~\ref{raisex1example} to $ y^2 - Q(x) $,
where $ p \ne 2 $, $ Q(x) $ in $ R[x] $ is
of degree $ 2g+1 $, and its reduction in $ k[x] $ has degree
$ 2g+1 $ and no multiple roots. (In other words, if we take $ R = W(k) $,
the Witt vectors of $ k $, then we are in the
situation studied in~\cite{Ked01}.)
Inverting $ 2 y $, we obtain an open part corresponding to
$ R[x,y,z]/(y^2-Q(x), 2yz-1) $.  
We have 
$ \Jac = \begin{pmatrix}
 -Q'(x) & 2y &  0 \\
  0     & 2z & 2y \\
\end{pmatrix}$,
so we can take
$ P = \begin{pmatrix}
             0      & 0   \\
    z    & 0   \\
   -2z^3 & z \\
\end{pmatrix}$,
$ \Delta_2 $ the zero matrix, and
$ \Delta_3 = 
\begin{pmatrix}
  1    & 0  \\
 -2z^2 & 1  \\
\end{pmatrix}
$.

In order to find a lift $ \phi $ of the relative Frobenius on $ \ol A'$,
we have to find the solution
$ \sol = (s_2 , s_3 ) $ in $(\pi R\la x,y,z \ra^\dagger)^2 $ of $G(\S) = 0$, where
\begin{alignat*}{1}
   G_2(\S)
& = 
   \bigl(y^{\pa} + z^{\pa} S_2 \bigr)^2 - Q^\s(x^{\pa}) - \(y^2 - Q(x)\)^{\pa} - \(2yz-1\)^{\pa} S_2
\\
\intertext{and}
   G_3(\S)
& =
   2 \bigl(y^{\pa} + z^{\pa} S_2 \bigr) \bigl( z^{\pa}- 2 z^{3\pa}S_2+z^{\pa} S_3 \bigr)-1
              - \bigl( 2yz-1 \bigr)^{\pa} \bigl(1 - 2 z^{2\pa}S_2 + S_3\bigr)
\,.
\end{alignat*}
If $ \sol = (s_2,s_3) $ is the unique solution in $ ( \pi R\la x, y, z \ra ^\dagger )^2  $,
then $ y^{\pa} + z^{\pa} s_2 $ is the unique solution for $ \phi(y) $ 
in $ \Adag = R\la x, y, z \ra^\dagger / (f, 2yz-1 ) $ 
of $ \phi(y)^2 - Q^\s(x^{\pa}) = 0 $ that is
congruent to $ y^{\pa} $ modulo $ \pi $, hence it must coincide
with the explicit formula given in \cite{Ked01} when $ \pa = p $
and $ R = W(k) $.
\end{example}

\section{Expansions at the ends} \label{local-frobenius}

We now return to our curve as described in Section~\ref{sec:intro}.
We extend the base field $ K $ to $ \tK $, and $ R $ to $ \tR $,
so that $ C_{\tR} \setminus X_{\tR} $ is the union of the `missing
points' $ Q_i $.
Let $ \ee $ be one the ends of the rigid analytic space corresponding
to $ C_{\tK} \setminus X_{\tK} $ described in Section~\ref{sec:cup-res},
and fix $ Q $, one of the missing points that lies in the corresponding
residue disc~$ \D $.

Because for this $ \ee $ we only need this section $ Q $,
we do not have to extend $ R $ to $ \tR $; it suffices to
replace $ R $ with a suitable $ R_Q \subseteq \tR $.
As this makes no difference to the proofs
we avoid this cumbersome notation and write $ R $ instead of
$ R_Q $ or $ \tR $.

In Theorem~\ref{globallift} we have constructed a lift $ \phi $
of $ \phibar $.  In order to calculate the contribution of $ \ee $
to the right hand side of~\eqref{localcup},
we could calculate the expansion of $ \phi(\eta) $ as follows.
We first apply $ \s $ to the coefficients of $ \eta $ in order,
compute the $ \phi(x_i) $ as well as their expanions, and substitute
the latter into $ \eta^\s $.
Instead, we never compute the $ \phi(x_i) $
but expand the $ x_i $ in the defining equations~\eqref{globalG}
and solve those.  This way we can obtain expansions of the $ \phi(x_i) $
directly, without the need of substituting expansions into expansions.
Another advantage is that we work with expressions that
contain only the local parameter, not all the variables $ x_i $.
Also, in practice the local expansions can converge on a larger annulus
than one might expect from the behaviour of the global $ \phi(x_i) $
(see Examples~\ref{runningexample3} and~\ref{runningexample4}).

The drawback is of course that we have to do solve the equations
for all ends $ \ee $, so if there are many of those, it may be
better to compute the $ \phi(x_i) $ globally first and then substitute
expansions of the $ x_i $.  In order to maintain this flexibility,
we also discuss how the estimates
on the coefficients in the global $ \phi(x_i) $ translate into
estimates on the coefficients in their local expansions.

Let $ t $ be a local equation of $ Q $ on $ C $ (as scheme, not
rigid analytic space), which we shall
also view as a parameter on $ \D $ and $ \ee $, and use it to
make the restriction map of rigid analytic
functions $ U_r $ (for a suitable $ r<1 $) to $ \ee $ 
explicit on $ \Adag $.

Let $ \OCq $ be the local ring at the reduction $ q $ of $ Q $
of $ C $.
Note that $ \OCq/(t) \iso R $,
and the completion of $ \OCq $ with respect to $ (t) $ is
isomorphic with $ R[[t]] $.  We shall refer to the resulting
map $ \OCq \to R[[t]] $,
or any of the analogues described below, as the expansion map
at $ Q $.
If $ a $ in $ A $ is such that in $ \D $ it
only has a pole along $ Q $, then $ t^{-\ord_Q(a)} a $ is in $ \OCq $
and $ a $ has an expansion in $ \Rt $.

More generally, let
\begin{equation*}
   \Rthat
=
 \biggl\{ \sum_{m\in\Z} a_m t^m \textup{ with all } a_m \textup{ in } R \textup{ and } \lim_{m\to-\infty}a_m=0 \biggr\}
\end{equation*} 
be the $ \pi $-adic completion of $ \Rt $.  Then 
any element in $ \Rt $
that is not in $ \pi \Rt $ is in $ \Rthat^* $: we can write it
as $ t^d f - \pi g $ with $ f $ in $ R[[t]]^* $ and $ g $ in $ \Rt $,
which has inverse $ t^{-d} f^{-1} ( 1 + \sum_{m \ge 1} (\pi t^{-d} f^{-1} g)^m ) $.
If $ a $ is any element in $ A $, then using local equations
in $ \OCq $ of irreducible divisors on $ C $ containing $ q $,
one sees that there is some $ h $ in $ \OCq $ such that $ h a $
is in $ \OCq $.  Because $ a $ does not restrict to 0 on $ C_k $,
we can assume the same about $ h $.  Because the composition
$ \OCq \to R[[t]] \to k[[t]] $ descends to the expansion map on
$ \O_{C_k, q} $, it follows that $ h $ maps to a unit $ u $ in $ \Rthat $.
Then $ a $ has the expansion in $ \Rthat $ obtained by expanding $ h a $ and
multiplying by $ u^{-1} $.

As $ \Rthat $ is $ \pi $-adically complete, the expansion map 
$ A = R[\x]/(f_2,\dots,f_n) \to \Rthat $
at $ Q $ extends to $ \Rx / (f_2,\dots,f_n) \to \Rthat $.
We shall see later that this extension restricted
to $ \Adag $ takes values in a suitable subring
$ \Rtstar $ of $ \Rthat $.
The extension also induces a map
\begin{equation} \label{expansionmap}
 \ex : \Rx \to \Rthat
\,.
\end{equation}
We shall abuse notation and denote
by $ \xi $ the map
to $ \Rthat $ from any of $ \Rx $, $ \Rxdag $, and $ \Adag $.

In order to describe the image of $ \Rxdag $ under this
map, together with estimates, below, we introduce a subring $ \Rtstar $
of $ \Rthat $. 
We shall show in Proposition~\ref{expest} 
that $ \xi $ maps $ \Rxdag $
into $ \Rtstar $, together with a description for bounds on the
coefficients involved.

In many applications the expansions of the $ x_j $ will be 
in $ \Rt $. We therefore include statements
that deal with this case specifically,
namely Remarks~\ref{cruderemark} and~\ref{subtleremark}.

In order to describe our estimates on coefficients we introduce
the following subsets of $ \Rthat $.
Note that each element in them is a rigid function 
as described in~\eqref{eq:annulus-fun} on (a possibly
narrower) $ \ee $.

\begin{notation}
For any rational numbers $\a$ and $\b$ with $\a>0$ we let
\begin{equation*}
 \Rtab {\a} {\b}
=
\biggl\{\sum_{m\in\Z}{a_m t^m} \text{ with all } a_m \text{ in }R \text{ and } v(a_m)\ge  -\a m+\b \biggr\}
\,.
\end{equation*}
We also let
$ \Rtstar = \bigcup_{\a,\b} \Rtab {\a} {\b} $.
\end{notation}
It can be helpful to visualize the conditions on the $ a_m $ by drawing the
region in the plane in which the pairs $ (m,v(a_m)) $ for non-zero $ a_m $
can lie, as in Figure~\ref{Rabfigure}.

\medskip

\begin{figure}[htbp]
\begin{center}
\input{Rab.pstex_t}
\caption{}
\label{Rabfigure}
\end{center}
\end{figure}

The following is easily established.

\begin{lemma} \label{rtabprop}
  The subsets above have the following properties.
  \begin{enumerate}
  \item The elements in $ \Rtab {\a} {\b} $ converge for $ p^{-\a} <
    |t| < 1 $.

  \item $ \Rtab {\a_1} {\b_1} \times \Rtab {\a_2} {\b_2} \to \Rtab
    {\min(\a_1,\a_2)} {\b_1+\b_2} $ under multiplication in $ \Rthat
    $.\label{rtabprop2}

  \item $ \Rtstar $ is a subring of $ \Rthat $, as are the $\Rtab {\a}
    0 $.

  \item $ \Rtab {\a} 0 $ is $ \pi $-adically complete.

  \item The units of $ \Rtab {\a} 0 $ are those $ \sum_m a_m t^m $
    with $ a_0 $ in $ R^* $. (Write such an element as $ u+w $ with $
    u = \sum_{n\ge0} a_m t^m $ and $ w = \sum_{m<0} a_m t^m $. Then
    its inverse is $ u^{-1}(1 - u^{-1}w + u^{-2}w^2 - \cdots) $.)

  \item If $ c $ is an element of $ R $ with $ v(c) = \cc > 0 $, then $
    c \Rtab {\a} {\b} $ is contained in $ \Rtab {\a} 0 $ if $ \b + \cc
    \ge 0 $, and in $ \Rtab {-\cc\a/\b} 0 $ if $ \b + \cc < 0 $ 
   $($see Figure~\ref{Rabcfigure}\/$)$. Note that $ -\cc\a/\b < \a $ when $ \b +
    \cc < 0 $.
  \end{enumerate}
\end{lemma}

\begin{figure}[htbp]
\begin{center}
\input{Rabc.pstex_t}
\caption{}
\label{Rabcfigure}
\end{center}
\end{figure}

\begin{remark} \label{minplus}
For a finite subset $ T $ of $ \R \cup \{\infty\} $,
we define $ \min^+(T) $  as $ \min(T \cap \R_{>0}) $;
i.e., we ignore all negative numbers as well as $ \infty $.
Then this last property states that
$ c \Rtab {\a} {\b} $ is contained in $ \Rtab {\a'} 0 $ with
$ \a' = \min^+\{ \a, -\cc\a/\b \} $.
\end{remark}

We now fulfill an earlier promise, and show that
the expansion map maps $ \Rxdag $ to $ \Rtstar $.
In particular, each element in $ \Adag $ is mapped to a rigid
function on (a possibly narrower) $ \ee $.

\begin{prop} \label{expest}
The expansion map $ \ex $ in~\eqref{expansionmap} maps $ \Rxdag $ to $ \Rtstar $.
More precisely, 
if $ g(\x) = \sum_I a_I \x^I $
in $ \Rxdag $ with $ v(a_I) \ge D \cdot I + \d $
for some $ D $ in $ (\ratpos)^m $ and $ \d $ in $ \Q $, then
$ \ex(g(\x)) $ is in $ \Rtab {\cc} {\d} $ where $ \cc $ is obtained
as follows:
\begin{enumerate}
\item
if all $ \ex(x_i) $ are in $ \Rt $, take
$ \cc $ in $ \ratpos $ with $ d_i \ge - \cc \ord_t(x_i) $ for all $ i $;

\item
if each $ \ex(x_i) $ is in some $ \Rtab {\a} { \b_i} $, 
let $ \cc' $ in $ \ratpos $ be such that 
$ - d_i \le \b_i \a^{-1} \cc' $ for all $ i $, and take
$ \cc= \min\{\cc', \a\} $.
\end{enumerate}
\end{prop}

\begin{proof}
(1)
Write $ D = (d_1,\dots,d_n) $, and
let $ D' = (d_1',\dots,d_n') $ with $ d_i' = \ord_t(x_i) $.
Note that each $ \ex(x^I) $ is in $ t^{D' \cdot I} R[[t]] $.
we have $ v(a_I) \ge D \cdot I + \d \ge - \cc D' \cdot I + \d  $,
so that $ a_I t^{D' \cdot I} $ is in $ \Rtab {\cc} {\d} $.
The same then holds for each $ \ex(a_I x^I ) $, hence for $ \ex(g(x)) $.
(Of course, if all $ d_i' \ge 0 $ then we can take  $ \cc $ arbitrarily
large and recover that $ \ex(g(x)) $ is in $ R[[t]] $.)

(2)
Let $ B=(\b_1,\dots,\b_n) $ so that the expansion of $ a_I x^I $ is
in $ a_I \Rtab {\a} {B\cdot I} $
by  Lemma~\ref{rtabprop}\eqref{rtabprop2}.
Then the vertex for the corresponding region as in Figure~\ref{Rabcfigure} occurs at
$ (\a^{-1} B \cdot I, v(a_I)) $.
But this point is above $ (\a^{-1} B \cdot I, \max\{0, D \cdot I + \d \} ) $,
which is to the right of $ ( -\cc'^{-1} D \cdot I , \max\{0, D \cdot I + \d \} ) $.
As $ I $ varies, those last points all lie in the region for $ \Rtab {\cc'} {\d} $.
If $ \a \ge \cc' $ then the same holds for the regions for 
all $ a_I \Rtab {\a}, {B \cdot I} $.  If $ \a \le \cc' $ then it holds
if we enlarge our region to that of $ \Rtab {\a} {\d} $.
\end{proof}

\begin{remark}
Note that if at least one $ d_i'<0 $ in part~(1) of Proposion~\ref{expest},
then $ \cc = \min^+\{ - d_i/d_i' \} $ is in $ \ratpos $, and is
the best possible choice. If $ \cc = -d_j/d_j' $, then $ D \cdot I  = -\cc D' \cdot I $
for all $ I $  having non-zero $ j $-th entry and zeroes elsewhere.
So for this $ \cc $ the statement of part~(1) appears to be
optimal.

The same cannot be said for part~(2) in general, because the
estimate is based on the vertex at the bend in Figure~\ref{Rabcfigure},
which may not correspond to an actual point $ (m, v(a_m)) $ for
an element of $ \Rtstar $.
For example, suppose $ n = 1 $, $ D= (d_1) $, and take $ \a = 2 $, $ \b = -1 $.
Then Proposition~\ref{expest} gives us $ \cc' = 2d_1 $ as largest
possible $ \cc' $, and the result $ \ex(g(x)) $ is in $ \Rtab {\cc} {\d} $
with $ \cc = \min\{2 d_1, 1\} $. On the other hand, if $ e=1 $,
then $ \Rtab 2 -1 \subset \Rtab 1 0 $. Taking $ g(x) $ in this
larger set, we can now take any $ \cc' $ in $ \ratpos $.
The result $ \ex(g(x)) $) lies in $ \Rtab 1 {\d} $.
(Note that $ d_1 $ drops out in this example because $ \Rtab 1 0 $
is a ring, so all $ \ex(x)^i $ are in it.)
\end{remark}

\begin{example} \label{runningexample3}
Let us return to the estimates obtained in Example~\ref{runningexample2}.
There we had an element $ \sum_{i,j} b_{i,j} x^i y^j $ with
$ v(b_I) \ge \frac{i+j}{6p} + \frac12 $,
so that $ D = (\frac{1}{6p} , \frac{1}{6p} ) $ and $ \d = \frac12 $.
Both $ x $ and $ y $ have poles of order 1 at each of the two
points at infinity, so the largest $ \cc $ we can take is $ \frac{1}{6p} $.
Then Proposition~\ref{expest}(1) states that $ \ex(\sum_{i,j} b_{i,j} x^i y^j ) $
is in $ \Rtab \frac{1}{6p} \frac12 $.
\end{example}

Note that we could compute the expansions of the $ \phi(x_i) $
constructed in Theorem~\ref{globallift0},
by first computing the $ \ex(x_j) $ and substituting those into
the $ \phi(x_i) $.  However, unless there are many ends,
it should be much more efficient
if we can compute the expansions of the
$ \phi(x_i) $ directly from their definition.
That this can be done is the content of Theorem~\ref{locallift0}.
In Theorem~\ref{locallift} we shall discuss estimates on the coefficients in the expansions
obtained by this method.
Note that the global estimates obtained in Theorem~\ref{globallift}
give us estimates on the expansions as well by applying
Proposition~\ref{expest}, but the two estimates can be quite different
(see Examples~\ref{runningexample3} and~\ref{runningexample4}).

\begin{theorem} \label{locallift0}
Let $ P $ and $ \Delta^2,\dots,\Delta^n $ and
$ G(\S) $ be as in Theorem~\ref{globallift},
and let $ \phi $ be the resulting $ \s $-linear endomorphism
of $ \Adag $.
Then the expansions at $ Q $ of the $ \phi(x_i) $ can be computed directly.
More precisely, if $ H(\S) $ is obtained from~\eqref{globalG}
by applying $ \ex $ to the coefficients, then there is a unique
solution $ \soll $ in $ (\pi \Rthat)^{n-r} $ of $ H(\S) = 0 $,
and $ \ex(\phi(\x)) = \ex(\psi(\x)) + \ex(\psi(P)) \soll $.
\end{theorem}

\begin{proof}
Recall that
$ \phi $ is induced by the $ \s $-linear endomorphism
of $ \Rxdag $ mapping $ g(\x) $ to $ g^\s(\psi(\x) + \psi(P) \sol ) $,
with $ \sol $ the unique solution in $(\pi \Rxdag)^{n-r}$ of
$ G(\S) = 0 $. 
So $ \ex(\phi(\x)) = \ex(\psi(\x)) + \ex(\psi(P)) \ex(\sol) $.
Because $ \ex(\sol) $ is in $ (\pi \Rtstar )^{n-r} $
by Proposition~\ref{expest}, it suffices to show that  
$ H(\S) = 0 $ has a unique solution $ \soll $ in
$ (\Rthat)^{n-r} $.

By that proposition the coefficients in $ H(\S) $
are in $ \Rtstar $, and $ H(\S) $
has inherited the following properties from $ G(\S) $:
\begin{itemize}
\item
$ H(0) \equiv 0 $ modulo $ \pi \Rtstar $;

\item
$ \Jac_H(0)  = \ex\( \Jac_G(0) \)$,
hence $\Jac_H(0) \equiv\Id_{n-r} $  modulo $ \pi \Rtstar $.
\end{itemize}
Applying Hensel's lemma for the $ \pi $-adically complete ring
$ \Rthat $ finishes the proof.
\end{proof}

\begin{remark}\label{locnewtonrem}
Note that applying $ \ex $ to the coefficients of~\eqref{globalG}
kills the terms involving the $ f_j $.
In particular, $ H(\S) $ is obtained by applying $ \ex $ to the
coefficients in $ f^\s (\psi(x)+\psi(P)\S) $, hence is determined by $ f $ and
$ P $.
The solution $ \soll $ we then obtain as
the appropriate solution of $ f^\s (\ex(\psi(x)) + \ex(\psi(P))\S ) = 0 $.
\end{remark}

In order to give estimates on the coefficients involved in the
solution $ \soll $ described in Theorem~\ref{locallift0}, we need
some lemmas and remarks. The reader should think of those as
the local analogue of Lemma~\ref{polytopelemma}.

\begin{lemma} \label{hensel1}
Let $ G_{r+1}(\S), \dots, G_n(\S) $ in $ \Rtab {\a} 0 [\S] $ for some $ \a>0 $
be of total maximal degree $ N > 0 $ in the variables $ \S = (S_{r+1},\dots,S_n) $.
Assume that each $ G_j(0) $ is in $ \pi^b \Rtab {\a} 0 $ for
some integer $ b \ge 1 $, and
that the determinant of $ \Jac_G(0) $ is in $ \Rtab {\a} 0 ^* $.
Then there is a unique solution $ \sol $ in $ ( \pi^b \Rtab {\a} 0 )^{n-r} $.
\end{lemma}

\begin{proof}
Apply Hensel's lemma to the equation $ G(\S) = 0 $, starting
with $ z_0 = (0,\dots,0) $ as approximate solution, and observe that
under Newton iteration $ z_{i+1} = z_i - \Jac_G(z_i)^{-1} G(z_i) $ we
stay in $ \pi^b \Rtab {\a} 0 $ all the time.
\end{proof}

\begin{remark}
Note that this is sharp for polynomials
of the form $ (S+1)^{N-1} (S-z) $ with $ z $ in $ \pi \Rtab {\a} 0 $.
\end{remark}

\begin{lemma} \label{hensel2}
Let $ H_{r+1}(\S),\dots,H_n(\S) $ be in $ \Rtstar [\S] $
of total maximal degree $ N \ge 1 $. Assume that there 
exist $ \a_l>0 $ and $ \b_l $ for $ l = 0, \dots, N $, and an integer
$ a\ge1 $, such that
\begin{itemize}
\item
all  $ H_j(0) $ are in $ \pi^a \Rtab {\a_0} {\b_0} $;

\item
the entries of $ \Jac_{H}(0) $ are in $ \Rtab {\a_1} 0 $
and its determinant is in $ \Rtab {\a_1} 0 ^* $;

\item
the homogeneous parts of degree $ l $ of all $ H_j(\S) $ 
are in $ \Rtab {\a_l} {\b_l} [\S]$ for $ l = 2, \dots, N $.
\end{itemize}
Then $ H(\S) = 0 $ has a unique solution $ \sol $ with coordinates
in $ \pi \Rtstar $.

In fact, if we write such a coordinate as
$ \sum_m a_m t^m $, then we have the following bound.
For any $ \nu $ in $ \ratpos $ satisfying $ \nu < a $, let
\begin{equation} \label{alphaprime}
\a_\nu' = \min{}\!^+\left\{
   \a_0,\dots,\a_N, 
    - \frac{(a - \nu ) \a_0}{e\b_0} , 
    - \frac{\nu \a_2}{e\b_2},
    - \frac{2 \nu \a_3}{e\b_3},
   \dots,
    - \frac{(N-1)\nu \a_N}{e\b_N}
\right\}
\,.
\end{equation}
Then $ v(a_m) \ge \max\{ 0, -\a_\nu' m \} + \frac{\nu}{e}
$. 
\end{lemma}

\begin{proof}
That there is a unique solution with coordinates in $ \pi \Rthat $ is again a consequence
of Hensel's lemma, since all $ \Rtab {\a} {\b} $ are in the $ \pi $-adically
complete ring $ \Rthat $.

Now fix a $ \nu $ in $ \ratpos $ with $ \nu < a $, and 
let $ \e $ in $ \ratpos $ satisfy $ \e < a - \nu $.
We shall be using Convention~\ref{ratconv} again.
Let $ G_j(\S') = \pi^{-\nu} H_j(\pi^\nu \S') $ for $j=r+1,\ldots,n$.
Then $ G_j(0) $ is in $ \pi^{a-\nu} \Rptab {\a_0} {\b_0} $,
$ \Jac_G(0) = \Jac_H(0)$ has entries in $ \Rptab {\a_1} 0 $ and determinant
in $ \Rptab {\a_1} 0 ^* $, and
the homogeneous parts of degree $ l $ of all the $ G_j(\S') $
are in $ \pi^{(l-1)\nu} \Rptab {\a_l} {\b_l} [\S'] $ for $ l=2,\dots,N $.
We can apply Lemma~\ref{hensel1}
(but with $ R $ replaced with $ R' $, $ \pi^a $ with $ \pi^\e $,
and $ \a $ with $ \a' $), provided that 
$ \pi^{a - \nu - \e } \Rptab {\a_0} {\b_0} $,
$ \Rptab {\a_1} 0  $, and the
$ \pi^{(l-1)\nu} \Rptab {\a_l} {\b_l} $ for $ l=2,\dots,N $ are
all in $ \Rptab {\a'} 0 $.
By Remark~\ref{minplus}, this is the case
when the following hold simultaneously:
\begin{itemize}
\item
$ \a' \le \min^+\{ \a_0 , - (a - \nu - \e ) \a_0/(e\b_0) \} $;

\item
$ \a' \le \a_1 $;

\item
$ \a' \le \min^+\{ \a_l , - (l-1)\nu \a_l/(e\b_l) \} $
for $ l=2,\dots,N $.
\end{itemize}
Therefore we can certainly take
\begin{equation*}
\a' = \a_{\nu,\e}' = \min{}\!^+\left\{
   \a_0,\dots,\a_N, 
    - \frac{(a - \nu - \e ) \a_0}{e\b_0} , 
    - \frac{\nu \a_2}{e\b_2},
   \dots,
    - \frac{(N-1)\nu \a_N}{e\b_N}
\right\}
\,.
\end{equation*}

From Lemma~\ref{hensel1} we obtain that the vector equation $ G(\S') = 0 $
has a unique solution $ \sol' $ with coordinates in $ \pi^\e \Rptab {\a'} 0 $.
Then $ \pi^{\nu} \sol' $ is a solution of $ H(\S) = 0 $ in
$ \pi^{\e+\nu} R'_*((t)) \subset \pi^a R'_*((t))$. But from Hensel's lemma 
in $ R'_*((t)) $, we see that $ H(\S) = 0 $
has a unique solution with coordinates in $ \pi^a \widehat{R'((t))} $,
hence $ \sol = \pi^\nu \sol' $ has coordinates in
$ \pi^{\nu+\e} \Rptab {\a'} 0 $.
As the coordinates are actually in $R$, 
letting $ \e $ go to zero finishes the proof.
\end{proof}

We can now obtain our main estimates for the expansions of the
$ \phi(x_i) $.

\begin{theorem} \label{locallift}
Let $ \soll $ be the unique solution of $ H(\S) = 0 $ constructed
in Theorem~\ref{locallift0}.
Then we have the following estimates for the coefficients in
$ \soll $.

\begin{enumerate}
\item
If there are positive integers $ a $ and $ b $ such that
the coefficients of the homogeneous parts of degree $ l $
are in $ \pi^a t^{-d_0} R[[t]] $ for $ l=0 $,
in $ R+ \pi^b t^{-d_1} R[[t]] $ for $ l=1 $,
in $ t^{-d_l} R[[t]] $ for $ l = 2, \dots, N $,
and the determinant of $ \Jac_H(0) $ is in $ R^* +  \pi^b t^{-d_1} R[[t]] $,
then for every $ \nu $ in $ \ratpos $ with $ \nu < a $, 
we have that $ \soll $ is in $ \Rtab {\a_\nu} {\frac{\nu}{e}} $,
where
\begin{equation*}
\a_\nu = 
\frac1e
  \min{}\!^+\left\{
      \frac{b}{ d_1} ,
      \frac{(a - \nu )}{ d_0} , 
      \frac{\nu}{ d_2},
      \frac{2 \nu }{ d_3},
   \dots,
      \frac{(N-1)\nu }{ d_N}
\right\}
\,.
\end{equation*}

\item
If in the statement of (1) we replace $ t^{-d_l} R((t)) $  with
$ \Rtab {\a} {- d_l \a} $ for $ l=0,\dots,N $, then the coordinates
of $ \soll $ are in $ \Rtab { \min\{\a, \a_\nu \} } {\frac{\nu}{e}} $,
with $ \a_\nu $ as in (1), again for all $ \nu $ in $ \ratpos $ with
$ \nu < a $.
\end{enumerate}
\end{theorem}

\begin{proof}
Clearly, $ t^{-d_l} R[[t]] \subset \Rtab {\a_l} {\b_l} $
if $ - d_l \ge \frac{\b_l}{\a_l} $. We take $ \b_l = d_l \a_l $
with $ \a_l $ very large for $ l \ne 1 $.
For $ l=1 $ we note that
$ \pi^b t^{-d_1} R[[t]] \subset \Rtab {\min^+\{\frac{b}{e d_1}\}} 0 $
(which we interpret as $ R[[t]] $ if $ d_1 \le 0 $).
We now apply Lemma~\ref{hensel2} to $ H(\S) $ as in Theorem~\ref{locallift0}.
Then~\eqref{alphaprime} simplifies
to the given expression for $ \a_\nu $ (except if all $ d_l \le 0 $,
in which case $ \a = + \infty $; but~\eqref{alphaprime}
can be made arbitrarily large in the same way). This proves part~(1).

For part~(2), we observe that
$ \pi^b \Rtab {\a} {-d_1 \a} \subset \Rtab { \min^+\{ \a , \frac{b}{e d_1}\} } 0 $
and apply Lemma~\ref{hensel2} to $ H(\S) $ as in Theorem~\ref{locallift0}.
\end{proof}

\begin{remark}  \label{cruderemark}
In explicit examples one can try to maximize the bound given
in Theorem~\ref{locallift}, but as a crude estimate, let us 
assume $ N \ge 2 $, and take $ \cc > 0 $ and $ \d \ge 0 $ such that
$ d_l \le (l-1)\cc + \d  $ for $ l=2,\dots,N $.
Then
\begin{equation*}
\frac1e
  \min{}\!^+\left\{
      \frac{b}{ d_1} ,
      \frac{(a - \nu )}{ d_0} , 
      \frac{\nu}{\cc + \d}
\right\}
\le \a_\nu
\end{equation*}
because $ \frac{(l-1) \nu }{ (l-1) \cc + \d} $  increases with
$ l $.
If $ d_0 \le 0 $ then we can let $ \nu $ approach $ a $ and obtain
$ \tilde \a= \frac1e \min{}\!^+\{ \frac{b}{ d_1} , \frac{a}{\cc + \d} \} $,
so that $ \soll $ has coordinates in
$ \Rtab {\a_\nu} {\frac{\nu}e } \subseteq \Rtab {\tilde\a} {\frac{a}e} $.
If $ d_0 > 0 $, then we equate the last two entries and solve for $ \nu $,
which gives $ \tilde\nu = \frac{a (\cc+\d)}{d_0+\cc+\d} $.
With $ \tilde\a = \frac1e \min{}\!^+\{ \frac{b}{ d_1} , \frac{a}{d_0 +\cc + \d} \} $
we have $ \a_{\tilde\nu} \ge \tilde\a $, and $ \soll $ has coordinates
in $ \Rtab {\tilde\a} {\frac{\tilde\nu}e} $.
\end{remark}

\begin{remark} \label{subtleremark}
Although Theorem~\ref{locallift} and Remark~\ref{cruderemark} give estimates when all
coordinates in Lemma~\ref{hensel2} have entries in $ \Rt $, they
do not take into account the coefficients involved in $ R $.
One can sometimes obtain better estimates by following the method
of the proof of Theorem~\ref{locallift}
and applying Lemma~\ref{hensel1} directly.  Namely,
in case~(1) of the theorem,
take $ \nu  $ in $ \ratpos $ with $ \nu < a $ and 
consider $ G(\S') = \pi^{-\nu} H(\pi^\nu \S) $.
We can then determine an $ \a_\nu $ as in the lemma by taking
the minimum of $ -v(a_i)/i $ over all $a_i t^i $
with $ i<0 $ in all coefficients of $ G(\S') $.
(Note that we replace $ \pi^a $ with $ \pi^\e $ for some $ \e $
in $ \ratpos $ again and let $ \e $ approach 0.)
It follows that the solution $ \tsoll $ has coordinates in
$ \pi^{-\nu} \Rtab {\a_\nu} 0 $.
Varying $ \nu $ we can select an $ \a_\nu $
that is optimal, or close to optimal.
Similar considerations apply in case~(2) of the theorem.
\end{remark}

\begin{example} \label{runningexample4}
Let us return to Example~\ref{runningexample1} and obtain local
estimates.
In Example~\ref{runningexample3}
we derived local esttimates from the global one in
Example~\ref{runningexample2}.
Recall that $ p \ne 2 $.
There are two points at infinity, $ [1,1,0] $ and $ [1,-1,0] $. 
With local parameter $ t = 1/x $ we find the expansions $ x(t) = t^{-1} $
and $ y(t) = \pm t^{-1} \sqrt{1- t^2} = \pm t^{-1} (1 - \frac12 t^2 - \frac18 t^4 -\frac{1}{16} t^6 - \cdots  ) $.
Using those in~\eqref{runexeq} we obtain
\begin{equation*}
   H(S)
 =  
   4^{-1} A(t) S^2 + A(t) S  - A(t) + 1
\end{equation*}
with $ A(t) = x(t)^{2p} - y(t)^{2p} = x(t)^{2p} - (x(t)^2-1)^p = p t^{-2p+2} + \cdots $
in $ 1 + p t^{-2p+2} \Rt $.
Then $ \a_\nu $ in Theorem~\ref{locallift}(1)
becomes $ \min\left\{\frac{1}{2p-2}, \frac{1-\nu}{2p-2}, \frac{\nu}{2p-2}\right\} $.
The best possible is when $ \nu = \frac12 $, 
so that the solution $ \ts $ that we want lies in
$ \Rtab {\frac{1}{4p-4}} {\frac12} $. The final estimate of Remark~\ref{cruderemark}
also gives this if we take $ \d = 2p-2 $ and let $ \cc $ approach~0.
But if we let $ m > 0 $, extend $ R $ to $ R' $ by working inside
a totally ramified of degree $ m $ over $ \Q_p $, then we can
apply Lemma~\ref{hensel1} directly with $ \a = \frac{1-1/m}{2p-2} $ and $ b=1 $.
Letting $ m $ go to infinity we find that $ \ts $ is in $ \Rtab {\frac{1}{2p-2}} {0} $.
Multiplying with $ 2 ^{-1} x_i(t)^p $ or $ 2^{-1} y_i(t)^p $ we obtain the local
expansions of $ \ex(\sum_{i,j} b_{i,j} x^i y^j ) $ of Example~\ref{runningexample3}
again.  Remembering that every coefficient in $ s $ contains
a factor $ p $, we find that the local expansion is in $ \Rtab {\frac{1}{2p-2}} {\frac{-p}{2p-2}} $.
This compares quite favourably with the estimates 
obtained in Example~\ref{runningexample2}, which were derived
directly from estimates on the global Frobenius $ \phi $.
\end{example}

\section{Examples of the local Frobenius} \label{local-examples}

In this section we revisit some of the examples in Section~\ref{global-examples}.
In particular, we investigate the case of an hyperelliptic curve as
in Example~\ref{kedlayaexample}, leaving out either the point
at infinity, or all the Weierstrass points.
We work out those cases mostly as an illustration of the differences
between leaving out as few points as possible, or opting for
localizing but imposing $ \phi(x) = x^{\pa} $.

The reader should bear in mind that for those curves,
using the closed formula as in \cite{Ked01}
for $ \phi $ with $ \phi(x) = x^{\pa} $, one can certainly get more precise information
about the expansions $ \ex(\phi(y)) $ than by our general methods.
Also, due to the low degree in $ y $ of the defining equation,
the problem of computing $ \phi(y) $ or its expansion
for this $ \phi(x) $
is of a rather different nature than in the case of a more general
curve.

Let the notation and assumptions be as in Example~\ref{kedlayaexample}.

\begin{example}
Let $ X $ to be the open affine corresponding to $ R[x,y]/(f) $,
so that we leave out only the point at infinity.
Since $ \ol Q $ has no multiple roots, there exist polynomials $ \ol{a}(x) $ of
degree at most $ 2g $ and $ \ol{b}(x) $ of degree at most $ 2g-1 $
in $ k[x] $ such that $ -\ol{a}(x) \ol{Q'}(x) + \ol b(x) \ol Q(x) = 1  $ in $ k[x] $.
Then
\begin{equation*}
	\ol a~\ol{f_x} + 2^{-1}y\ol b~\ol{f_y} = 1 + \ol b~\ol f
\end{equation*}
in $ k[x,y] $.
We lift $ \ol{a} $ and $ \ol{b} $ to $ a $ and $ b $ of degree
at most $ 2g $ and $ 2g-1 $ in $ R[x] $, so that
$ a $ and $ 2^{-1} y b $ have a pole at infinity of order at most $ 4g $ and $ 6g-1 $ respectively.

Using those, the $ H(S) $ defined in Theorem~\ref{locallift0} becomes
\begin{alignat*}{1}
	H(S) &= \bigl(y(t)^{\pa} + \tfrac{1}{2} y(t)^{\pa} b^\s (x(t)^{\pa} ) S \bigr)^2 -
			Q^{\s}\bigl(x(t)^{\pa} + a^\s(x(t)^{\pa})S\bigr)\\
         &= \sum_{l=0}^{2g+1} H_l S^l\,.
\end{alignat*}
Choosing a parameter $ t $ centred at the missing point,
the expansions $ x(t) $ and $ y(t) $ are in $ t^{-2}R[[t]] $ and 
$ t^{-2g-1}R[[t]] $ respectively. 
So $ H_0 $ is in $ \pi t^{-2\pa(2g+1)}R[[t]] $,
$ H_1 $ is in
$ R^* + \pi t^{-8\pa g}R[[t]] $, and 
$ H_l $ is in 
$ t^{-2\pa(l(2g-1) + 2g+1)}R[[t]] = t^{-2\pa((l-1)(2g-1)+4g)}R[[t]] $
for $ l = 2,\cdots,2g+1 $.

Applying Theorem~\ref{locallift0}, there is a unique solution $ \ts $ in $
\pi\Rthat $ of $ H(S) = 0 $. Moreover, Remark~\ref{cruderemark} gives us the 
following estimate on $ \ts $. For every 
$ \nu $ in $ \ratpos $, with $ \nu < 1 $, $ \ts $ is in 
$ \Rtab {\a_\nu} {\frac{\nu}{e}} $, where
\begin{equation*}
\frac1e
  \min{}\!^+\left\{
  \frac{1}{8\pa g},
  \frac{(1 - \nu )}{ 2\pa(2g+1)} , 
  \frac{\nu}{2\pa(6g-1)}
\right\}
\le \a_\nu \,.
\end{equation*}
Equating the last two entries and solving for $ \nu $
gives $ \tilde\nu = \frac{6g-1}{8g} $.
With $ \tilde\a = \frac1e \frac{1}{16\pa g} $
we have $ \a_{\tilde\nu} \ge \tilde\a $, and $ \ts $ is in 
$ \Rtab {\tilde\a} {\frac{\tilde\nu}e} $.

Note that $ x(t)^{\pa} $ is in 
$ \Rtab {\tilde\a} {-2\pa\tilde\a} $, hence
$ a^\s(x(t)^{\pa}) $ is in
$ \Rtab {\tilde\a_\nu} {-4\pa g\tilde\a} $,
and $ \ex(\phi(x)) = x(t)^{\pa} + a^\s(x(t)^{\pa}) \ts $ is in 
$ \Rtab {\tilde\a_\nu} {\b} $ with 
$ \b = \min\{-2\pa\tilde\a, -4\pa g \tilde\a + \frac{\tilde\nu}e \} $.
Similarly, $ \ex(\phi(y)) = y(t)^{\pa} + \frac{1}{2}b^\s(x(t)^{\pa}) \ts $ is
in $ \Rtab {\tilde\a} {\b'} $ with
$ \b' = \min\{-\pa ( 2g+1) \tilde\a, -2\pa (2g-1)\tilde\a + \frac{\tilde\nu}e \} $.
\end{example}

\begin{example}
Let us now invert $ 2y $ as in Example~\ref{kedlayaexample},
so that we work
with the open affine corresponding to $ R[x,y,z]/(y^2 - Q(x), 2yz - 1) $
and the missing points are the $ 2g+2 $ Weierstrass points.
The vector $ H(\S) $ defined in Theorem~\ref{locallift0} has entries
\begin{alignat*}{1}
	H_2(\S) & = \bigl(y(t)^{\pa} + z(t)^{\pa} S_2\bigr)^2 - Q^\s(x(t)^{\pa}) 
				= H_{2,0} + H_{2,1} S_2 + H_{2,2} S_2^2 \\
\intertext{and}
	H_3(\S) & = 2 \bigl(y(t)^{\pa} + z(t)^{\pa} S_2\bigr) \bigl( z(t)^{\pa}-
				2z(t)^{3\pa}S_2+z(t)^{\pa} S_3 \bigr)-1
\,.
\end{alignat*}

As in Example~\ref{raisex1example}, 
the first condition involves only $ S_2 $,
and by Lemma~\ref{hensel2}, $ H_2(S_2) = 0 $ has a unique solution $ \ts_2 $
in $ \pi\Rtstar $. To give estimates on $ \ts_2 $, we need to study three
distinct cases.

Case 1: the missing point is the point at infinity. Then the expansions 
$ \ex(x) $, $ \ex(y) $ and $ \ex(z) $ are in $ t^{-2}R[[t]] $,
$ t^{-2g-1}R[[t]] $ and $ t^{2g+1}R[[t]] $ respectively. Hence
$ H_{2,0} $ is in $ \pi t^{-2\pa(2g+1)}R[[t]] $, $ H_{2,1} $
is in $ R^* $, and $ H_{2,2} $ is in $ t^{2\pa(2g+1)} R[[t]] $. 
Letting $ \nu $ approach $ 0 $ and taking $ b $ large in Theorem~\ref{locallift}
applied to $ H_2 $, we get that $ \ts_2 $ is in 
$ \Rtab {\a} 0 $ for $ \a = \frac{1}{2e\pa(2g+1)} $
(and of course all of its coefficients are in $ \pi R $).
Then $ \ex(\phi(y)) = y(t)^{\pa} + z(t)^{\pa} \ts_2 $ is in 
$ \Rtab {\a} {-\frac 1 {2e}} $.
We also have 
\begin{equation} \label{phizcalc}
 \ex(\phi(z)) = \frac12 \ex( \phi(y))^{-1} = 
\frac12 y(t)^{-\pa} \frac{1}{1 + 2^{\pa} z(t)^{2\pa} \ts_2}
\,,
\end{equation}
which is in $ t^{\pa(2g+1)} \Rtab {\a} 0 $
because $ z(t)^{2\pa} \ts_2 $ is in $ \pi \Rtab {\a} 0 $.

Case 2: the missing point is a Weierstrass point $ (a,0) $.
Here we can choose $ y $ to be the local parameter, and the expansions $
\ex(x) $ and $ \ex(z) $ are in $ R[[t]] $ and $ t^{-1}R[[t]] $ respectively.
Therefore $ H_{2,0} $ is in $ \pi R[[t]] $, $ H_{2,1} $ is
in $ R^* $, and $ H_{2,2} $ is in $ t^{-2\pa}R[[t]] $.
Now letting $ \nu $ approach $ a = 1 $ and taking $ b $ large
in Theorem~\ref{locallift} applied to $ H_2 $, we get that $ \ts_2 $ is in 
$ \Rtab {\a} {\frac{1}{e}} $ with $ \a = \frac 1 {2e\pa} $.
Then $ \ex(\phi(y)) $ is in $ t^{-\pa} \Rtab {\a} {\frac1e} $
as this contains both $ y(t)^{\pa} $ and $ z(t)^{\pa} \ts_2 $.
Computing $ \ex( \phi(z) ) $ as in~\eqref{phizcalc}
we see that it is in $ t^{-\pa} \Rtab {\a} 0 $.
\end{example}

\section{Finite precision estimates} \label{sec:estimates}

In this section we explain how to use the methods that were discussed in
Section~\ref{local-frobenius} to get the cup products required
in Method~\ref{cupmethod} up to a given precision.
Here knowing $ c $ in $ K $ or $ \tK $ up to precision $ N $
means that we have an explicit
$ \tilde c $ in $ K $ or $ \tK $ with $ v(c - \tilde c) \ge N $.
In order to simplify notation, for $ a $ in $ \Q $, we shall write $ \Ia a $
for $ \{ x\text{ in } K \text{ with } v(x)\ge a \}$, so that
we want to find a representative $ \tilde c $ of a class in $ K / \Ia N  $.

We place ourselves in the situation of Section~\ref{local-frobenius},
so fix an end $ \ee $ and a local parameter $t$ for the corresponding
residue disc $ \D $.  As in that section, we write $ R $ for
what might be an extension of the original~$ R $.

We shall use the images of the $ \Rtab {\a} {\b} $
with finite precision for the coefficients.

\begin{notation}
For $ \a $, $ N $ in $ \ratpos $ and $ \b $ in $ \Q $, we define the set
\begin{equation}\label{stab}
  \Stab {\a} {\b} {N}
=
  \biggl\{\sum_m \ol{d_m} t^m \text{ with } \sum_m d_m t^m \text{ in } R_{\a,\b}((t)) \biggr\}
\subseteq
  S^{N}((t))
\,,
\end{equation}
where we take the coefficients in the quotient ring $ S^{N} = R/ \Ia N $.
\end{notation}

As $ d_m $ satisfies $ v(d_m) \ge - m \a + \b $,
its image $ \ol{d_m} $ in $ S^{N} $ is trivial if
$ m \le (\b - N )/\a $.

Our basic computational problem is as follows.
Given forms $\eta$ and $\omega$ of
the second kind with expansions
\begin{equation*}
  \eta = \sump m {} a_m t^m \dd t, \quad \o= \sump m {} b_m t^m \dd t 
\end{equation*}
in the local parameter $t$, we need to determine
\begin{equation*}
  \res_\ee \o \int \eta = \sump m {} \frac{a_m b_{-m-2}}{m+1}
\,.
\end{equation*}

We want to show that this residue can be obtained  up to precision $N$
via a finite object with which we can compute. For this we shall
use the quotients
\begin{equation*}
 \Stab {\a} {\b} {N} / t^L =  \Stab {\a} {\b} {N} / t^L \Stab {\a} {\b} {N} 
\end{equation*}
for positive exponents $L$.  Note that we have products
\begin{equation} \label{Lmult}
   \Stab {\a} {\b_1} {N} / t^L
\times 
   \Stab {\a} {\b_2} {N} / t^L
\to
   \Stab {\a} {\b_1 + \b_2} {N} / t^L
\end{equation}
that are compatible with the products on the $\Rtab {\a} {\b_i} $.

Let us first analyse more closely the structure of
$ \Stab {\a} {\b} {N} / t^L $.

\begin{lemma}
If $\sum a_m t^m $ is in $ \Stab {\a} {\b} {N} / t^L $, then
we have $ a_m $ in $ \Ia {f(m)} / \Ia {g(m)} $, with
$ f(m)=\max(0,-\a m +  \b) $
and
$  g(m) = \max(f(m),\min(N,f(m-L)))$.
\end{lemma}

\begin{proof}
The statement means more precisely that we have a map onto the above
set which is compatible with the obvious map from $ \Rtab {\a} {\b} $
to $R$ which extracts the coefficient $a_m$. The condition with $f$
is obvious from the lower bounds on coefficients in $ \Rtab {\a} {\b} $.
The condition with $g$ comes from the fact that we are multiplying
by $t^L$ and quotienting out by the result.
In particular, the coefficient of $t^m$ in the resulting class must be taken
modulo the possible coefficients of $t^{m-L}$. Finally, $g(m)$ has to be at
least as large as $f(m)$. This is because even though the precision
is capped at $N$, if we know it is $0$ to a higher precision then it
is definitely known to this higher precision (and noting that $f(m)\le f(m-L)$
because $\a,L>0$).
\end{proof}

For our estimates we shall assume for simplicity that $\b_i\le 0$. This
occurs in practice and can be assumed by at worse replacing a positive
$\b_i$ with $0$.

\begin{prop}\label{estprop}
Given $ N $ in $ \ratpos $,
the map   $(\o,\eta)\mapsto  \res_\ee \o \int \eta + \Ia N $ factors via
$ \Stab {\a_1} {\b_1} {N_1} / t^{L_1} \cdot \dd t \times \Stab {\a_2} {\b_2} {N_2} / t^{L_2} \cdot \dd t$
for suitable $N_1, N_2$ in $ \ratpos $ and
positive integers $L_1 , L_2$.
\end{prop}

\begin{proof}
We observe that we have a well-defined multiplicaton map
\begin{equation*}
  \Ia {a_1} / \Ia {b_1} \times  \Ia {a_2} / \Ia {b_2} 
\to
  \Ia {a_1+a_2} / \Ia {\min(a_1+b_2,a_2+b_1)} 
\,.
\end{equation*}
From this, the above lemma and the definition of the residue, it is
clear that we can factor the residue as required if we have for all
$m\ne -1$ that
\begin{equation*}
  \min(g_1(m) + f_2(-m-2) ,  f_1(m)+g_2(-m-2) ) 
   \ge N +v(m+1)
\,.
\end{equation*}
To achieve this, we start by observing that the left hand side is
greater than or equal to $f_1(m) + f_2(-m-2) $,
which is independent of any choice of $ N_i $ and $ L_i $.
This is at least $\max(-\a_1 m+\b_1,-\a_2(-m-2)+\b_2)$ and thus,
for sufficiently large $|m|$, the above inequality certainly holds
independently of the choice of the $L_i$ and $N_i$.
We therefore need to 
find them so that the condition is satisfied for the finitely
many remaining~$ m $.

To this end, we may first guarantee the condition after taking
$N_i=\infty$ and finding appropriate $L_i$. Then $N_i$ can be taken
sufficiently large to make sure that the inequalities still hold. As
noted before,  $f_i(m)\le f_i(m-L_i)$.
Thus, for the remaining $m$'s our goal is to choose
$L_i$ so that
\begin{equation*}
  \min(f_1(m-L_1)+f_2(-m-2),f_1(m) +f_2(-m-2-L_2))\ge N+v(m+1)
\,.
\end{equation*}
Clearly, for each fixed $m$, this will be achieved for sufficiently
large $L_1$ and $L_2$.
\end{proof}

For computational purposes we provide a way of finding the relevant
constants.

\begin{prop}
  The following algorithm provides constants $L_i$ and $N_i$
  satisfying the conditions of Proposition~\ref{estprop}.
\begin{enumerate}
\item
Find integers
$M_+ \ge 0 $ and $M_- \le -2 $
with $-\a_1 m+\b_1-\log_p(|m+1|)\ge N$
for $m< M_-$ and $-\a_2(-m-2)+\b_2-\log_p(m+1)\ge N$
for $m> M_+$. Define $ \Mlog = \log_p (\max(M_+ + 1, - M_- - 1 ) ) $.

\item
Let $L_1$ and $L_2$ be positive integers satisfying
the following conditions.
\begin{enumerate}
\item \label{76a}
Let $m_0$ be the smallest integer not equal to $ -1 $ with $-\a_1 m_0+\beta_1\le 0$.
Let $L_2$ satisfy $-\a_2(-m_0-2-L_2)+\b_2 \ge N + \Mlog $.

\item \label{76b}
Let $m_1$ be the largest integer not equal to $ -1 $ with $-\a_2 (-m_2-2)+\beta_2\le 0$.
Let $L_1$ satisfy $-\a_1(m_1-L_1)+\b_1 \ge N + \Mlog $.

\item \label{76c}
If $\a_2>\a_1$, then
$-\a_1 M_-+\b_1 -\a_2 (-M_--2-L_2)+\b_2\ge N + \Mlog $.

\item \label{76d}
If $\a_1>\a_2$, then
$-\a_1(M_+-L_1)+\b_1-\a_2(-M_+-2)+\b_2\ge N + \Mlog $.

\end{enumerate}

\item \label{76e}
Take $N_1=N_2 \ge N + \Mlog $.
\end{enumerate}
\end{prop}

\begin{proof}
Note that the first step implies that $ f_1(m) + f_2(-m-2) \ge N + \log_p |m+1| $
for all $ m \ne -1 $ not in $ [M_- , M_+] $, as was done in the
proof of Proposition~\ref{estprop}.
For the remaining $ m \ne -1 $,
we may replace the term $v(m+1)$ in the required inequality by $ \Mlog $,
which is the maximum of $\log_p(|m+1|)$ for all such~$ m $. Our goal
is then to choose 
\begin{align}
  L_1 \text{ such that } f_1(m-L_1)+f_2(-m-2)&\ge N+ \Mlog, \label{l1cond}\\
  L_2 \text{ such that } f_1(m) +f_2(-m-2-L_2)&\ge N+ \Mlog  \label{l2cond}
\end{align}
for all $ m \ne -1 $ in $ [M_- , M_+] $.

First consider the smallest integer $m_0 \ne -1 $ for which
$f_1(m_0)=0$. The condition on $L_2$ coming from \eqref{l2cond} at $m_0$
is $f_2(-m_0-2-L_2) \ge N + \Mlog $, which follows from
\eqref{76a}. For $m>m_0$ we still have $f_1(m)=0$ while $f_2(-m-2-L_2)
$ is increasing in $m$, so the condition~\eqref{l2cond} continues to hold.

Now consider~\eqref{l2cond} for $m<m_0$. It is then implied by
\[
  -\a_1 m+\b_1 -\a_2 (-m-2-L_2)+\b_2\ge N + \Mlog
\,, 
\]
which we already imposed for $ m=m_0 $.
This is equivalent with
\[L_2\ge \a_2^{-1}( (\a_1-\a_2) m-\b_1 -2\a_2 -\b_2 + N + \Mlog ).\]
Suppose that
$\a_2 \le \a_1$. Then this induces a weaker bound on $L_2$ when $m$
decreases, so the existing condition coming from $m_0$ suffices.
On the other hand, for $\a_2>\a_1$ we should add an extra condition on $L_2$
coming from the smallest $m\ne -1$ in $ [M_-, M_+] $.
For this~\eqref{76c} suffices.

Similarly, we may consider the largest $m_1$
for which $f_2(-m_1-2)=0$. Then~\eqref{l1cond} for $m = m_1$ is
$f_1(m_1-L_1) \ge N + \Mlog $, which follows from \eqref{76b}, and it
continues to hold for all $m< m_1$.
For $m> m_1$, it is implied by
$-\a_1(m-L_1)+\b_1-\a_2(-m-2)+\b_2\ge N + \Mlog $.
As before, we get the extra condition \eqref{76d}.

The estimate on the $N_i$ is obvious (and probably not quite optimal).
\end{proof}

The computations in Section~\ref{local-frobenius} will give us the
$\a_i$ and $\b_i$ for $\eta$ and $\o$ from which we can compute the
parameters $L_i$ and $N_i$. Then, as the computation of the local
expansion of $\phi(\o)$ does not involve a loss in precision, the
residue calculation can be done using
the $\Stab {\a_i} {\b_i} {N_i} /t^{L_i} \cdot \dd t$.

\section{Algorithm and implementation} \label{algorithm}

In this section we describe the resulting algorithm for point counting
and give a rough estimate 
for its complexity. It is hard to give a very
precise bound because this could vary significantly among different
types of curves. We have also made various simplifying assumptions. We
shall be using the soft $O$ notation $\ot$, meaning that logrithmic factors
are neglected compared with polynomial ones, so that for example
$O(l\log^9(l)) = \ot(l)$.

We first recall some basic facts about the zeta function of $ C_k $
from~\cite{Weil49}.
As mentioned in Method~\ref{zeta-method}, the zeta-function of $ C_k $
is obtained as $ Z(T) = \frac{P_1(T)}{(1-T)(1-qT)} $, where
$ P_1(T) = a_0+\dots+a_{2g}T^{2g} $ in $ \Z[T] $ is $ \det(1-T M') $.
Then $ a_0=1 $ and $ a_{2g-i}=q^{g-i}a_i $ for $i=0,\dots,g$,
so that only need to know $ a_1, \dots, a_g $.
Moreover, if we write $P_1(T)=\prod_{j=1}^{2g} (1-\a_jT)$
in $ \mathbb{C}[T] $, then all $ |\a_j| = q^{\frac 1 2} $, therefore
we have $|a_i| \le {2g \choose i} q^{\frac{i}2} $ for $ i=1,\dots,g $.
Hence, it suffices to know $ \a_1,\dots,\a_{g} $
up to precision $ N = \log_p (2 {2g \choose g} q^{g/2}) $ in order to determine
their correct value in $ \Z $. Asymptotically we have $N=\ot(g^2 l)$.

To simplify matters, we shall assume that $K$ is unramified over $\Q_p$
and that $C\setminus X$ conists of a finite number of disjoint
$R$ sections. Not assuming
this probably does not change the complexity much because one is
typically working over a larger extension but at the same time the
results of the computation, being Galois conjugates of one another,
can be computed once for a bunch of points.

We let $q=p^l$ be the size of the residue field $ k $.
We are interested in asymptotics in $l$, so we shall
be assuming that $l$ is very large compared with $p$ and $g$. We shall
assume that other required data: Number of missing residue discs,
degrees of defining functions, degrees of functions showing up in the
matrix $P$, are linear in $g$. Note that there may well be situations
where this is over pessimistic. For example, for all hyperelliptic
curves with odd degree models we can manage with just the residue disc
at infinity. We
shall also not keep track on the dependence on the number of defining
equations, as this tends to be very small.

To compute the zeta function we need to compute the entries in $M$ to
precision $N$. Since we are assuming that $K$ is unramified we have at
our disposal the results of Berthelot  \cite[(2.1.4) of Chapter~VII]{Ber74},
to the effect that the cup
product pairing on crystalline cohomology is perfect. One further
knows that Frobenius
acts on integral crystalline cohomology and that $ \hcr^1(C_k/R) \cong
\hdr^1(C/R)$~\cite[3.4.2]{Ill74}. This implies that
the entries of both matrices $M_1$ and $M_2$ from Method~\ref{cupmethod} are
integral and the
determinant of $M_1$ is invertible, provided we start with a basis for
the integral de Rham cohomology of $C$. Thus, both matrices are still
required at precision $N$. 
We ignore here the issue of finding an integral basis, but this
is in practice easily done using expansions of polar parts.

Using the contents of Sections~\ref{sec:cup-res}, \ref{global-frobenius},
\ref{local-frobenius}, and~\ref{sec:estimates}, we can now give
an algorithm that computes the numerator $ P_1(T) $ of the zeta
function of a curve.

\medskip

\begin{algorithm}
\ 
\smallskip

INPUT:

\begin{itemize}
\item 
A presentation of an $ R $-algebra $ A = R[\x]/(f) $ that satisfies
Assumption~\ref{Aassumption}, and such that $ A $ corresponds to an open affine 
$ X \subset C $ with $ C \setminus X $ consisting of the union of disjoint sections $ Q_i $.

\item
Matrices $ P $ and $ \Delta^{r+1},\dots,\Delta^n $ in $ \M^{n,n-r}(R[\x]) $ and
$ \M^{n-r}(R[\x]) $ respectively, such that
$ \Jac_{(f)}\times P \equiv \textup{Id}_{n-r} + \sum_{j=r+1}^n f_j\Delta^j $
modulo $ p $.

\item 
For every missing point $ Q $, the local expansions $ \ex(x_i) $ at $ Q $.

\item
A set of representatives $ \o_1,\dots,\o_{2g} $ in $ \Of {\Adag/R} $,
for a basis of the image of $ \hcr^1(C_k/R) $ inside 
$ H_\rig^1 (\X/K) $.
\end{itemize}

\medskip

{\bf Step 1: preliminary precision estimates.}
\begin{enumerate}
\item
Determine the required precision 
$ N = \log_p (2 {2g \choose g} q^{g/2}) $. 

\item
For all missing points $Q=Q_1,\dots,Q_r$ do

\begin{enumerate}
\item
Compute, in $ \Rthat[Z] $, $ H(\S) $ as defined in Theorem~\ref{locallift0}.

\item
Using the estimates from Theorem~\ref{locallift} and Remark~\ref{cruderemark},
determine $ \cc $ in $ \ratpos $ and $ \d $ in $ \Q $ such that the components of $ \soll $ lie in 
$ \Rtab {\cc} {\d} $.

\item
Using the equality 
$ \ex(\phi(x_i)) = x_i(t)^p + \sum_{j=r+1}^n \ex(\psi(P_{i,j})) \tilde s_i $, determine 
$ \tilde\a $ in $ \ratpos $ and $ \tilde \b_i $ in $ \Q $ such that each
$ \phi(x_i) $ has its expansion in $ \Rtab {\tilde\a} {\tilde\b_i} $.

\item
For every form $ \o $, compute $ \a $ in 
$ \ratpos $ and $ \b $ in $ \Q $, such that $ \phi(\o) $ has its
expansion in $ \Rtab {\a} {\b} $.

\item
Using those, compute, as explained in
Section~\ref{sec:estimates}, the precision $L_i$ in $ t $ required for the various
residue computations.

\item
If the precision in $ t $ of the $ x_i(t) $ is not big enough, then exit with an
error, otherwise continue.

\item
Also compute the biggest $N_i$ needed in the computations of the
expansions of the $ \phi(\o) $.
\end{enumerate}

\end{enumerate}

\medskip

{\bf Step 2: computation of the matrix $ M_1 $.}

\medskip

{\bf Step 3: computation of the local lifts of $ \phi $.}

\smallskip
For all missing points $Q=Q_1$,\dots,$Q_r$ do
\begin{enumerate}
\item
Use Newton iteration to compute, up to precision $ N_i $,
the solution $ \soll $ of $ H(\S) = 0 $ with $ \soll\equiv 0 $ modulo $ p $.

\item
Deduce from that the Laurent series
$ \ex(\phi(x_i)) = x_i(t)^p + \sum_{j=r+1}^n \ex(\psi(P_{i,j})) \tilde s_j $ 
up to the same precision.
\end{enumerate}

\medskip

{\bf Step 4: computation of the matrix $ M_2 $.}

\smallskip
For all forms $\o$, $\eta$ in the given basis and for all missing point $Q$ do
\begin{enumerate}
\item
With notation as in Section~\ref{sec:estimates},
find the class of $ \eta $ in 
$ \Stab {\a_1} {\b_1} {N_1} / t^{L_1}  \cdot\dd t $,
represented in the form
$ \sum \ol{a}_m t^m\cdot\dd t $.

\item
Use the multiplication in~\eqref{Lmult}
to get the class of $ \phi(\o) $ in 
$ \Stab {\a} {\b} {N_2} / t^{L_2} \cdot\dd t $,
represented in the form
$ \sum  \ol{c}_m t^m\cdot\dd t $.

\item
Compute the value
$ \sum_m' \frac{\ol{a}_m \ol{c}_{-m-2}}{m+1} $ 
of $ \res_\ee \phi(\o) \int \eta $ in $ K/\Ia N $.

\item
Sum these residues over all $ \ee $ to obtain the entry in $ R / \Ia N $ of $M_2$
corresponding to $\langle\eta,F(\o)\rangle_U$.
\end{enumerate}

\medskip

{\bf Step 5: computation of the zeta function $ Z(T) $.}
\begin{enumerate}
\item
Compute the product $ M_1^{-1}M_2 $ corresponding to the matrix $ M $ of the
action of the $\s$-linear Frobenius, up to precision $ N $.

\item 
With $ q = p^r $, compute 
$ M' = \sigma^{r-1}(M)\times \sigma^{r-2}(M)\times \cdots \times \sigma(M)\times M $, 
the matrix of the action of the linear Frobenius, up to precision $ N $.

\item
Lift the coefficients $ \ol{a}_1,\dots, \ol{a}_g $ of the characteristic
polynomial 
$ \sum_{i=0}^{2g} \ol{a}_i T^i $ of $ M' $ to the unique integral numbers $ a_i $
satisfying
$ |a_i| \le {2g \choose i} q^{\frac i 2} $.

\item
Let $ a_0 = 1$ and compute $ a_{2g-i} $ as $ q^{g-i} a_i $ for $ i=0,\dots,g $.

\end{enumerate}

\medskip

OUTPUT:
If the starting precision is high enough, the numerator $ P_1(T) $ of the zeta function of $ C_k $.
\end{algorithm}

This algorithm has the following complexity.

\begin{prop}
The asymptotic complexity of this algorithm is $\ot(p l^3)$, where the
$\ot$ term depends polynomially on the genus.
\end{prop}

\begin{proof}
In this proof we shall also make an attempt to estimate the dependency
on the genus.
The computation of the matrices $M_1$ and $M_2$ involve a cup product
computation, which in turn decomposes into certain residue
computations as described in Section~\ref{sec:estimates}. We consider the
computations of $M_2$ as these are clearly more time consuming.

The computation is ``essentially'' done integrally. In other words,
considering the residue computation in Section~\ref{sec:estimates} the integration
introduces denominators, but these are
fairly mild.
For the asymptotics this introduces logarithmic factors that
will be swallowed by the $\ot$-notation.

By Section~\ref{sec:estimates} the complexity of the residue
computation is controlled by the
parameter $\a$ of overconvergence. Indeed, if our forms are in $  \Rtab {\a}
{0} $ (it is clear that from the point of view of the asymptotic
complexity the parameter $\b$ may be neglected), then all residue
computation may be done in the quotient rings $ \Stab {\a} {0} {N_i} /t^{L_i} $
of the rings defined in~\eqref{stab}, 
where $L_i$ is approximately $N/\a$. and $N_i$ is approximately $N$. Element in  $ \Stab {\a} {0} {N_i} /t^{L_i} $  are
Laurent series, truncated from both above and below at $L_i$, with
coefficients in $R/\Ia {N_i} $ (with some divisibility conditions for the
negative coefficients and modulo lower powers of $p$ for the positive
ones). The complexity of operations in this ring, including the final
residue operation, and using fast arithmetic, is $\ot(L_i)$ operations
in $R/ \Ia {N_i} $ which is $\ot(L_iN_i)$ operations in the residue field $k$. As
this has size $p^l$, operations take $\ot(\log(p^l))= \ot(l)$, taking
into account that $p$ will occur polynomially in the complexity. We can take $L_1= L_2=L$ and $N_i=N$ for evaluating the complexity.

Let us now count the number of ring operations required to compute the
elements of $M_2$. For each residue disc we first need to compute the
expansion of the $\phi(x_i)$. Here, we first need to compute the
coefficients for the required equations and then carry out Newton
iterations to solve them. As the convergence of the solution is with
respect to the $ p $-adic topology, the number of iterations is
proportional to the log of the $p$-adic precision, which is $N$, hence
ultimately to $\log(l)$. After that, we have to substitute the
expressions for $\phi(x_i)$ into the forms, an operation which has a
complexity of operations in $ \Stab {\a} {0} {N} /t^{L} $ proportional
to the total degree of the
defining expressions for these forms. This will be roughly
quadratic in $g$. By our asymptotic assumptions (here we are being
rather rough as we are assuming that $g$ is small compared with
$\log(l)$ and not just $l$), the dominant term will be the Newton
iteration. Each such iteration involves a computation, controlled by
the size of $f$, which is polynomial in $g$ (for a plane curve, the
total degree is of order $\sqrt{g}$, so the total number of
multiplications and additions required to carry out the Newton
iteration is of order $O(g\sqrt{g})$). Overall, the computation is done in
about $O(g^{\mu_1} \log(l))$ operations,
here for a plane curve $\mu_1=3/2$. This has to be further multiplied
by the number of residue discs, assumed to be  $O(g)$. 
Absorbing $\log(l)$ into the soft $O$, we get an overall complexity of
\begin{equation*}
  \ot(g^{\mu_1+1} LNl) = \ot(g^{\mu_1+1} l N^2/\a)
= \ot(g^{\mu_1+1} l N^2/\a)  = \ot(g^{\mu_1+5}l^3/\a)
\,.
\end{equation*}
Thus, the overall complexity depends on the size of $\alpha$. This can
be estimated using Part (1) of Theorem~\ref{locallift}. We assume that
the entries in the matrix $P$ will have poles of order $O(g^{\mu_2})$
at the removed points $Q_i$ (for plane curves one may take
$\mu_2 = 1/2$). The defining equation further involves
applying $\psi$ to the entries in $P$ (see Remark~\ref{locnewtonrem}),
multiplying the order of pole by $p$. Thus, overall we can expect $\a = 1/O(pg^{\mu_3})$,
which gives an overall complexity $ \ot(g^{\mu_1+\mu_3+5} l^3 p) $ for the
residue computation. Other required operations
fall within this bound~\cite{Ked01}.
\end{proof}


\end{document}